\documentclass{article}
\usepackage{amssymb}
\usepackage{amsfonts}
\usepackage{amsmath}

\input{tcilatex}
\begin{document}

\title{Intersection between the geometry of generalized Lie algebroids and
some aspects of interior and exterior differential systems}
\author{Constantin M. ARCU\c{S} \\
\ \ \\
[0pt] \ \ \\
[0pt] 
\begin{tabular}{c}
SECONDARY SCHOOL \textquotedblleft CORNELIUS RADU\textquotedblright , \\ 
RADINESTI VILLAGE, 217196, GORJ COUNTY, ROMANIA \\ 
e-mail: c\_arcus@yahoo.com, c\_arcus@radinesti.ro%
\end{tabular}%
}
\maketitle

\begin{abstract}
An exterior differential calculus in the general framework of generalized
Lie algebroids is presented. A theorem of Maurer-Cartan type is obtained.
All results with details proofs are presented and a new point of view over
exterior differential calculus for Lie algebroids is obtained. Using the
theory of linear connections of Ehresmann type presented in the first
reference, the identities of Cartan and Bianchi type are presented.
Supposing that any vector subbundle of the pull-back Lie algebroid of a
generalized Lie algebroid is \emph{interior differential system (IDS)} for
that generalized Lie algebroid, then the involutivity of the \emph{IDS} in a
theorem of Frobenius type is characterized. Extending the classical notion
of \emph{exterior differential system (EDS)} to generalized Lie algebroids,
then the involutivity of an \emph{IDS} in a theorem of Cartan type is
characterized. \ \ \ \ \bigskip\newline
\end{abstract}

\tableofcontents

\section{Introduction}

\ \ 

In general, if $\mathcal{C}$ is a category, then we denote $\left\vert 
\mathcal{C}\right\vert $ the class of objects and for any $A,B{\in }%
\left\vert \mathcal{C}\right\vert $, we denote $\mathcal{C}\left( A,B\right) 
$ the set of morphisms of $A$ source and $B$ target and $Iso_{\mathcal{C}%
}\left( A,B\right) $ the set of $\mathcal{C}$-isomorphisms of $A$ source and 
$B$ target. Let $\mathbf{Liealg},~\mathbf{Mod,}$ $\mathbf{Man}$ and $\mathbf{%
B}^{\mathbf{v}}$ be the category of Lie algebras, modules, manifolds and
vector bundles respectively.

We know that if 
\begin{equation*}
\left( E,\pi ,M\right) \in \left\vert \mathbf{B}^{\mathbf{v}}\right\vert ,
\end{equation*}
\begin{equation*}
\Gamma \left( E,\pi ,M\right) =\left\{ u\in \mathbf{Man}\left( M,E\right)
:u\circ \pi =Id_{M}\right\}
\end{equation*}
and 
\begin{equation*}
\mathcal{F}\left( M\right) =\mathbf{Man}\left( M,\mathbb{R}\right) ,
\end{equation*}
then $\left( \Gamma \left( E,\pi ,M\right) ,+,\cdot \right) $ is a $\mathcal{%
F}\left( M\right) $-module.

If \ $\left( \varphi ,\varphi _{0}\right) \in \mathbf{B}^{\mathbf{v}}\left(
\left( E,\pi ,M\right) ,\left( E^{\prime },\pi ^{\prime },M^{\prime }\right)
\right) $ such that $\varphi _{0}\in Iso_{\mathbf{Man}}\left( M,M^{\prime
}\right) ,$ then, using the operation 
\begin{equation*}
\begin{array}{ccc}
\mathcal{F}\left( M\right) \times \Gamma \left( E^{\prime },\pi ^{\prime
},M^{\prime }\right) & ^{\underrightarrow{~\ \ \cdot ~\ \ }} & \Gamma \left(
E^{\prime },\pi ^{\prime },M^{\prime }\right) \\ 
\left( f,u^{\prime }\right) & \longmapsto & f\circ \varphi _{0}^{-1}\cdot
u^{\prime }%
\end{array}%
\end{equation*}%
it results that $\left( \Gamma \left( E^{\prime },\pi ^{\prime },M^{\prime
}\right) ,+,\cdot \right) $ is a $\mathcal{F}\left( M\right) $-module and we
obtain the $\mathbf{Mod}$-morphism%
\begin{equation*}
\begin{array}{ccc}
\Gamma \left( E,\pi ,M\right) & ^{\underrightarrow{~\ \ \Gamma \left(
\varphi ,\varphi _{0}\right) ~\ \ }} & \Gamma \left( E^{\prime },\pi
^{\prime },M^{\prime }\right) \\ 
u & \longmapsto & \Gamma \left( \varphi ,\varphi _{0}\right) u%
\end{array}%
\end{equation*}%
defined by 
\begin{equation*}
\begin{array}{c}
\Gamma \left( \varphi ,\varphi _{0}\right) u\left( y\right) =\varphi \left(
u_{\varphi _{0}^{-1}\left( y\right) }\right) =\left( \varphi \circ u\circ
\varphi _{0}^{-1}\right) \left( y\right) ,%
\end{array}%
\end{equation*}%
for any $y\in M^{\prime }.$

If $M,N\in \left\vert \mathbf{Man}\right\vert ,$ $h\in Iso_{\mathbf{Man}%
}\left( M,N\right) $, $\eta \in Iso_{\mathbf{Man}}\left( N,M\right) $ and $%
\left( F,\nu ,N\right) \in \left\vert \mathbf{B}^{\mathbf{v}}\right\vert $
so that there exists 
\begin{equation*}
\begin{array}{c}
\left( \rho ,\eta \right) \in \mathbf{B}^{\mathbf{v}}\left( \left( F,\nu
,N\right) ,\left( TM,\tau _{M},M\right) \right)%
\end{array}%
\end{equation*}%
and an operation 
\begin{equation*}
\begin{array}{ccc}
\Gamma \left( F,\nu ,N\right) \times \Gamma \left( F,\nu ,N\right) & ^{%
\underrightarrow{\left[ ,\right] _{F,h}}} & \Gamma \left( F,\nu ,N\right) \\ 
\left( u,v\right) & \longmapsto & \left[ u,v\right] _{F,h}%
\end{array}%
\end{equation*}%
with the following properties:\bigskip

\noindent $\qquad GLA_{1}$. \emph{the equality holds good }%
\begin{equation*}
\begin{array}{c}
\left[ u,f\cdot v\right] _{F,h}=f\left[ u,v\right] _{F,h}+\Gamma \left(
Th\circ \rho ,h\circ \eta \right) \left( u\right) f\cdot v,%
\end{array}%
\end{equation*}%
\qquad \quad\ \ \emph{for all }$u,v\in \Gamma \left( F,\nu ,N\right) $\emph{%
\ and} $f\in \mathcal{F}\left( N\right) .$

\medskip $GLA_{2}$. \emph{the }$4$\emph{-tuple} $\left( \Gamma \left( F,\nu
,N\right) ,+,\cdot ,\left[ ,\right] _{F,h}\right) $ \emph{is a Lie} $%
\mathcal{F}\left( N\right) $\emph{-algebra,}

$GLA_{3}$. \emph{the }$\mathbf{Mod}$\emph{-morphism }$\Gamma \left( Th\circ
\rho ,h\circ \eta \right) $\emph{\ is a }$\mathbf{LieAlg}$\emph{-morphism of 
}%
\begin{equation*}
\left( \Gamma \left( F,\nu ,N\right) ,+,\cdot ,\left[ ,\right] _{F,h}\right)
\end{equation*}%
\emph{\ source and }%
\begin{equation*}
\left( \Gamma \left( TN,\tau _{N},N\right) ,+,\cdot ,\left[ ,\right]
_{TN}\right)
\end{equation*}%
\emph{target, }then the triple 
\begin{equation*}
\begin{array}{c}
\left( \left( F,\nu ,N\right) ,\left[ ,\right] _{F,h},\left( \rho ,\eta
\right) \right)%
\end{array}%
\leqno(1.1)
\end{equation*}%
is an object of the category $\mathbf{GLA}$\emph{\ }of generalized Lie
algebroids. (see $\left[ 1\right] $)

Using the exterior differential algebra of the generalized Lie algebroid $%
\left( 1.1\right) $ we develop the theory of Lie derivative, interior
product and exterior differential operator in Section $2.$ A new theorem of
Maurer-Cartan type is presented in this general framework. As any Lie
algebroid can be regarded as a generalized Lie algebroid, then a new point
of view over exterior differential calculus for Lie algebroids is obtained.
(see also: $\left[ 4,5,10,11\right] $) All the results are given with
detailed proofs. Using the $\left( \rho ,h\right) $-torsion and $\left( \rho
,h\right) $-curvature presented in $\left[ 1\right] ,$ we obtain the $\left(
\rho ,h\right) $-torsion and $\left( \rho ,h\right) $-curvature forms and
identities of Cartan and Bianchi type in Section $3.$

Using the \emph{Cartan's moving frame method}, there exists the following

\textbf{Theorem }(E. Cartan) \emph{If }$N\in \left\vert \mathbf{Man}%
_{n}\right\vert $\emph{\ is a Riemannian manifold and }$X_{\alpha
}=X_{\alpha }^{i}\frac{\partial }{\partial x^{i}},$\emph{\ }$\alpha \in 
\overline{1,n}$\emph{\ is an ortonormal moving frame, then there exists a
colection of }$1$\emph{-forms }$\Omega _{\beta }^{\alpha },~\alpha ,\beta
\in \overline{1,n}$\emph{\ uniquely defined by the requirements}%
\begin{equation*}
\Omega _{\beta }^{\alpha }=-\Omega _{\alpha }^{\beta }
\end{equation*}%
\emph{and }%
\begin{equation*}
d^{F}\Theta ^{\alpha }=\Omega _{\beta }^{\alpha }\wedge \Theta ^{\beta
},~\alpha \in \overline{1,n}
\end{equation*}%
\emph{where }$\left\{ \Theta ^{\alpha },\alpha \in \overline{1,n}\right\} $%
\emph{\ is the coframe.} (see $\left[ 12\right] ,$ p. $151$)

We know that an $r$\emph{-dimensional distribution on a manifold }$N$ is a
mapping $\mathcal{D}$ defined on $N,$ which assignes to each point $x$ of $N$
an $r$-dimensional linear subspace $\mathcal{D}_{x}$ of $T_{x}N.$ A vector
fields $X$ belongs to $\mathcal{D}$ if we have $X_{x}\in \mathcal{D}_{x}$
for each $x\in N.$ When this happens we write $X\in \Gamma \left( \mathcal{D}%
\right) .$

The distribution $\mathcal{D}$ on a manifold $N$ is said to be \emph{%
differentiable} if for any $x\in N$ there exists $r$ differentiable linearly
independent vector fields $X_{1},...,X_{r}\in \Gamma \left( \mathcal{D}%
\right) $ in a neighborhood of $x.$ The distribution $\mathcal{D}$ is said
to be \emph{involutive }if for all vector fields $X,Y\in \Gamma \left( 
\mathcal{D}\right) $ we have $\left[ X,Y\right] \in \Gamma \left( \mathcal{D}%
\right) .$

In the classical theory we have the following

\textbf{Theorem }(Frobenius)\textbf{\ }\emph{The distribution }$\mathcal{D}$%
\emph{\ is involutive if and only if for each }$x\in N$\emph{\ there exists
a neighborhood }$U$\emph{\ and }$n-r$\emph{\ linearly independent }$1$\emph{%
-forms }$\Theta ^{r+1},...,\Theta ^{n}$\emph{\ on }$U$\emph{\ which vanish
on }$\mathcal{D}$\emph{\ and satisfy the condition }%
\begin{equation*}
d^{F}\Theta ^{\alpha }=\Sigma _{\beta \in \overline{r+1,p}}\Omega _{\beta
}^{\alpha }\wedge \Theta ^{\beta },~\alpha \in \overline{r+1,n}.
\end{equation*}%
\emph{for suitable }$1$\emph{-forms }$\Omega _{\beta }^{\alpha },~\alpha
,\beta \in \overline{r+1,n}.$(see $\left[ 9\right] ,$ p. $58$)

Extending the notion of distribution, in Section $4,$ we obtain the
definition of an \emph{IDS} of a generalized Lie algebroid and a
characterization of the ivolutivity of an \emph{IDS} in a result of
Frobenius type is presented in \emph{Theorem 4.1. }In particular, $%
h=Id_{M}=\eta ,$ then we obtain the theorem of Frobenius type for Lie
algebroids. (see: $\left[ 2\right] ,$ p. $248$)

In Section $4$ of this paper we will show that there exists very close links
between \emph{EDS }and the geometry of generalized Lie algebroids. In the
classical sense, an \emph{EDS} is a pair $(M,\mathcal{I})$ consisting of a
smooth manifold $M$ and a homogeneous, differentially closed ideal $\mathcal{%
I}$ in the algebra of smooth differential forms on $M$. ( see: $[3,6,7,8])$
Extending the notion of \emph{EDS} to generalized Lie algebroids, the
involutivity of an \emph{IDS} in a result of Cartan type is presented in the 
\emph{Theorem 4.3. }In particular, $h=Id_{M}=\eta ,$ then we obtain the
theorem of Cartan type for Lie algebroids. (see: $\left[ 2\right] ,$ p. $249$%
)

Finally, in Section $5,$ we present a new direction by research in
Symplectic Geometry.

\bigskip

\section{Exterior differential calculus}

\ \ \ 

We propose an exterior differential calculus in the general framework of
generalized Lie algebroids. As any Lie algebroid can be regarded as a
generalized Lie algebroid, in particular, we obtain a new point of view over
the exterior differential calculus for Lie algebroids.

Let $\left( \left( F,\nu ,N\right) ,\left[ ,\right] _{F,h},\left( \rho ,\eta
\right) \right) \in \left\vert \mathbf{GLA}\right\vert $ be.

\textbf{Definition 2.1 }For any $q\in \mathbb{N}$ we denote by $\left(
\Sigma _{q},\circ \right) $ the permutations group of the set $\left\{
1,2,...,q\right\} .$

\textbf{Definition 2.2 }We denoted by $\Lambda ^{q}\left( F,\nu ,N\right) $
the set of $q$-linear applications 
\begin{equation*}
\begin{array}{ccc}
\Gamma \left( F,\nu ,N\right) ^{q} & ^{\underrightarrow{\ \ \omega \ \ }} & 
\mathcal{F}\left( N\right) \\ 
\left( z_{1},...,z_{q}\right) & \longmapsto & \omega \left(
z_{1},...,z_{q}\right)%
\end{array}%
\end{equation*}%
such that 
\begin{equation*}
\begin{array}{c}
\omega \left( z_{\sigma \left( 1\right) },...,z_{\sigma \left( q\right)
}\right) =sgn\left( \sigma \right) \cdot \omega \left( z_{1},...,z_{q}\right)%
\end{array}%
\end{equation*}%
for any $z_{1},...,z_{q}\in \Gamma \left( F,\nu ,N\right) $ and for any $%
\sigma \in \Sigma _{q}$.

The elements of $\Lambda ^{q}\left( F,\nu ,N\right) $ will be called \emph{%
differential forms of degree }$q$ or \emph{differential }$q$\emph{-forms}$.$

\textit{\noindent Remark 2.1}\textbf{\ }If $\omega \in \Lambda ^{q}\left(
F,\nu ,N\right) $, then $\omega \left( z_{1},...,z,...,z,...z_{q}\right) =0.$
Therefore, if $\omega \in \Lambda ^{q}\left( F,\nu ,N\right) $, then $\omega
\left( z_{1},...,z_{i},...,z_{j},...z_{q}\right) =-\omega \left(
z_{1},...,z_{j},...,z_{i},...z_{q}\right) .$

\textbf{Theorem 2.1 }\emph{If }$q\in N$\emph{, then }$\left( \Lambda
^{q}\left( F,\nu ,N\right) ,+,\cdot \right) $\emph{\ is a} $\mathcal{F}%
\left( N\right) $\emph{-module.}

\textbf{Definition 2.3 }If $\omega \in \Lambda ^{q}\left( F,\nu ,N\right) $
and $\theta \in \Lambda ^{r}\left( F,\nu ,N\right) $, then the $\left(
q+r\right) $-form $\omega \wedge \theta $ defined by%
\begin{equation*}
\begin{array}{ll}
\omega \wedge \theta \left( z_{1},...,z_{q+r}\right) & =\underset{\sigma
\left( q+1\right) <...<\sigma \left( q+r\right) }{\underset{\sigma \left(
1\right) <...<\sigma \left( q\right) }{\tsum }}sgn\left( \sigma \right)
\omega \left( z_{\sigma \left( 1\right) },...,z_{\sigma \left( q\right)
}\right) \theta \left( z_{\sigma \left( q+1\right) },...,z_{\sigma \left(
q+r\right) }\right) \vspace{3mm} \\ 
& \displaystyle=\frac{1}{q!r!}\underset{\sigma \in \Sigma _{q+r}}{\tsum }%
sgn\left( \sigma \right) \omega \left( z_{\sigma \left( 1\right)
},...,z_{\sigma \left( q\right) }\right) \theta \left( z_{\sigma \left(
q+1\right) },...,z_{\sigma \left( q+r\right) }\right) ,%
\end{array}%
\end{equation*}%
for any $z_{1},...,z_{q+r}\in \Gamma \left( F,\nu ,N\right) ,$ will be
called \emph{the exterior product of the forms }$\omega $ \emph{and}~$\theta
.$

Using the previous definition, we obtain

\textbf{Theorem 2.2 }\emph{The following affirmations hold good:} \medskip

\noindent 1. \emph{If }$\omega \in \Lambda ^{q}\left( F,\nu ,N\right) $\emph{%
\ and }$\theta \in \Lambda ^{r}\left( F,\nu ,N\right) $\emph{, then}%
\begin{equation*}
\begin{array}{c}
\omega \wedge \theta =\left( -1\right) ^{q\cdot r}\theta \wedge \omega .%
\end{array}%
\leqno(2.1)
\end{equation*}

\noindent 2. \emph{For any }$\omega \in \Lambda ^{q}\left( F,\nu ,N\right) $%
\emph{, }$\theta \in \Lambda ^{r}\left( F,\nu ,N\right) $\emph{\ and }$\eta
\in \Lambda ^{s}\left( F,\nu ,N\right) $\emph{\ we obtain}%
\begin{equation*}
\begin{array}{c}
\left( \omega \wedge \theta \right) \wedge \eta =\omega \wedge \left( \theta
\wedge \eta \right) .%
\end{array}%
\leqno(2.2)
\end{equation*}

\noindent 3. \emph{For any }$\omega ,\theta \in \Lambda ^{q}\left( F,\nu
,N\right) $\emph{\ and }$\eta \in \Lambda ^{s}\left( F,\nu ,N\right) $\emph{%
\ we obtain }%
\begin{equation*}
\begin{array}{c}
\left( \omega +\theta \right) \wedge \eta =\omega \wedge \eta +\theta \wedge
\eta .%
\end{array}%
\leqno(2.3)
\end{equation*}

\noindent 4. \emph{For any }$\omega \in \Lambda ^{q}\left( F,\nu ,N\right) $%
\emph{\ and }$\theta ,\eta \in \Lambda ^{s}\left( F,\nu ,N\right) $\emph{\
we obtain }%
\begin{equation*}
\begin{array}{c}
\omega \wedge \left( \theta +\eta \right) =\omega \wedge \theta +\omega
\wedge \eta .%
\end{array}%
\leqno(2.4)
\end{equation*}

\noindent 5. \emph{For any} $f\in \mathcal{F}\left( N\right) $, $\omega \in
\Lambda ^{q}\left( F,\nu ,N\right) $ \emph{and }$\theta \in \Lambda
^{s}\left( F,\nu ,N\right) $\emph{\ we obtain }%
\begin{equation*}
\begin{array}{c}
\left( f\cdot \omega \right) \wedge \theta =f\cdot \left( \omega \wedge
\theta \right) =\omega \wedge \left( f\cdot \theta \right) .%
\end{array}%
\leqno(2.5)
\end{equation*}%
\noindent

\textbf{Theorem 2.3 }\emph{If }%
\begin{equation*}
\Lambda \left( F,\nu ,N\right) =\underset{q\geq 0}{\oplus }\Lambda
^{q}\left( F,\nu ,N\right) ,
\end{equation*}%
\emph{then }$\left( \Lambda \left( F,\nu ,N\right) ,+,\cdot ,\wedge \right) $%
\emph{\ is a} $\mathcal{F}\left( N\right) $\emph{-algebra.} This algebra
will be called \emph{the exterior differential algebra of the vector bundle }%
$\left( F,\nu ,N\right) .$

\emph{Remark 2.2 }If $\left\{ t^{\alpha },~\alpha \in \overline{1,p}\right\} 
$ is the coframe associated to the frame $\left\{ t_{\alpha },~\alpha \in 
\overline{1,p}\right\} $ of the vector bundle $\left( F,\nu ,N\right) $ in
the vector local $\left( n+p\right) $-chart $U$, then 
\begin{equation*}
\begin{array}{c}
t^{\alpha _{1}}\wedge ...\wedge t^{\alpha _{q}}\left( z_{1}^{\alpha
}t_{\alpha },...,z_{q}^{\alpha }t_{\alpha }\right) =\frac{1}{q!}\det
\left\Vert 
\begin{array}{ccc}
z_{1}^{\alpha _{1}} & ... & z_{1}^{\alpha _{q}} \\ 
... & ... & ... \\ 
z_{q}^{\alpha _{1}} & ... & z_{q}^{\alpha _{q}}%
\end{array}%
\right\Vert ,%
\end{array}%
\leqno(2.6)
\end{equation*}%
for any $q\in \overline{1,p}.$

\emph{Remark 2.3 }If $\left\{ t^{\alpha },~\alpha \in \overline{1,p}\right\} 
$ is the coframe associated to the frame $\left\{ t_{\alpha },~\alpha \in 
\overline{1,p}\right\} $ of the vector bundle $\left( F,\nu ,N\right) $ in
the vector local $\left( n+p\right) $-chart $U$, then, for any $q\in 
\overline{1,p}$ we define $C_{p}^{q}$ exterior differential forms of the
type 
\begin{equation*}
\begin{array}{c}
t^{\alpha _{1}}\wedge ...\wedge t^{\alpha _{q}}%
\end{array}%
\end{equation*}%
such that $1\leq \alpha _{1}<...<\alpha _{q}\leq p.$

The set%
\begin{equation*}
\begin{array}{c}
\left\{ t^{\alpha _{1}}\wedge ...\wedge t^{\alpha _{q}},1\leq \alpha
_{1}<...<\alpha _{q}\leq p\right\}%
\end{array}%
\end{equation*}%
is a base for the $\mathcal{F}\left( N\right) $-module 
\begin{equation*}
\left( \Lambda ^{q}\left( F,\nu ,N\right) ,+,\cdot \right) .
\end{equation*}

Therefore, if $\omega \in \Lambda ^{q}\left( F,\nu ,N\right) $, then%
\begin{equation*}
\begin{array}{c}
\omega =\omega _{\alpha _{1}...\alpha _{q}}t^{\alpha _{1}}\wedge ...\wedge
t^{\alpha _{q}}.%
\end{array}%
\end{equation*}

In particular, if $\omega $ is an exterior differential $p$-form $\omega ,$
then we can written 
\begin{equation*}
\begin{array}{c}
\omega =a\cdot t^{1}\wedge ...\wedge t^{p},%
\end{array}%
\end{equation*}%
where $a\in \mathcal{F}\left( N\right) .$

\textbf{Definition 2.4 }If 
\begin{equation*}
\begin{array}{c}
\omega =\omega _{\alpha _{1}...\alpha _{q}}t^{\alpha _{1}}\wedge ...\wedge
t^{\alpha _{q}}\in \Lambda ^{q}\left( F,\nu ,N\right)%
\end{array}%
\end{equation*}%
such that 
\begin{equation*}
\begin{array}{c}
\omega _{\alpha _{1}...\alpha _{q}}\in C^{r}\left( N\right) ,%
\end{array}%
\end{equation*}%
for any $1\leq \alpha _{1}<...<\alpha _{q}\leq p$, then we will say that 
\emph{the }$q$\emph{-form }$\omega $\emph{\ is differentiable of }$C^{r}$%
\emph{-class.}

\textbf{Definition 2.5 }For any $z\in \Gamma \left( F,\nu ,N\right) $, the $%
\mathcal{F}\left( N\right) $-multilinear application 
\begin{equation*}
\begin{array}{c}
\begin{array}{rcl}
\Lambda \left( F,\nu ,N\right) & ^{\underrightarrow{~\ \ L_{z}~\ \ }} & 
\Lambda \left( F,\nu ,N\right)%
\end{array}%
,%
\end{array}%
\end{equation*}%
defined by%
\begin{equation*}
\begin{array}{c}
L_{z}\left( f\right) =\Gamma \left( Th\circ \rho ,h\circ \eta \right)
z\left( f\right) ,~\forall f\in \mathcal{F}\left( N\right)%
\end{array}%
\leqno(2.7)
\end{equation*}%
and 
\begin{equation*}
\begin{array}{cl}
L_{z}\omega \left( z_{1},...,z_{q}\right) & =\Gamma \left( Th\circ \rho
,h\circ \eta \right) z\left( \omega \left( z_{1},...,z_{q}\right) \right) \\ 
& -\overset{q}{\underset{i=1}{\tsum }}\omega \left( \left( z_{1},...,\left[
z,z_{i}\right] _{F,h},...,z_{q}\right) \right) ,%
\end{array}%
\leqno(2.8)
\end{equation*}%
for any $\omega \in \Lambda ^{q}\mathbf{\ }\left( F,\nu ,N\right) $ and $%
z_{1},...,z_{q}\in \Gamma \left( F,\nu ,N\right) ,$ will be called \emph{the
covariant Lie derivative with respect to the section }$z.$

In the particular case of Lie algebroids, $\left( \eta ,h\right) =\left(
Id_{M},Id_{M}\right) ,$ then \emph{the covariant Lie derivative with respect
to the section }$z$ is defined by 
\begin{equation*}
\begin{array}{c}
L_{z}\left( f\right) =\Gamma \left( \rho ,Id_{M}\right) z\left( f\right)
,~\forall f\in \mathcal{F}\left( M\right)%
\end{array}%
\leqno(2.7^{\prime })
\end{equation*}%
and 
\begin{equation*}
\begin{array}{cl}
L_{z}\omega \left( z_{1},...,z_{q}\right) & =\Gamma \left( \rho
,Id_{M}\right) z\left( \omega \left( z_{1},...,z_{q}\right) \right) \\ 
& -\overset{q}{\underset{i=1}{\tsum }}\omega \left( \left( z_{1},...,\left[
z,z_{i}\right] _{F},...,z_{q}\right) \right) ,%
\end{array}%
\leqno(2.8^{\prime })
\end{equation*}%
for any $\omega \in \Lambda ^{q}\mathbf{\ }\left( F,\nu ,M\right) $ and $%
z_{1},...,z_{q}\in \Gamma \left( F,\nu ,M\right) .$

In addition, if $\rho =Id_{TM},$ then \emph{the covariant Lie derivative
with respect to the vector field }$z$ is defined by 
\begin{equation*}
\begin{array}{c}
L_{z}\left( f\right) =z\left( f\right) ,~\forall f\in \mathcal{F}\left(
M\right)%
\end{array}%
\leqno(2.7^{\prime \prime })
\end{equation*}%
and 
\begin{equation*}
\begin{array}{cl}
L_{z}\omega \left( z_{1},...,z_{q}\right) & =z\left( \omega \left(
z_{1},...,z_{q}\right) \right) \\ 
& -\overset{q}{\underset{i=1}{\tsum }}\omega \left( \left( z_{1},...,\left[
z,z_{i}\right] _{TM},...,z_{q}\right) \right) ,%
\end{array}%
\leqno(2.8^{\prime \prime })
\end{equation*}%
for any $\omega \in \Lambda ^{q}\mathbf{\ }\left( TM,\nu ,M\right) $ and $%
z_{1},...,z_{q}\in \Gamma \left( TM,\nu ,M\right) .$

\textbf{Theorem 2.4 }\emph{If }$z\in \Gamma \left( F,\nu ,N\right) ,$ $%
\omega \in \Lambda ^{q}\left( F,\nu ,N\right) $\emph{\ and }$\theta \in
\Lambda ^{r}\left( F,\nu ,N\right) $\emph{, then}%
\begin{equation*}
\begin{array}{c}
L_{z}\left( \omega \wedge \theta \right) =L_{z}\omega \wedge \theta +\omega
\wedge L_{z}\theta .%
\end{array}%
\leqno(2.9)
\end{equation*}

\emph{Proof. }Let $z_{1},...,z_{q+r}\in \Gamma \left( F,\nu ,N\right) $ be
arbitrary. Since 
\begin{equation*}
\begin{array}{l}
L_{z}\left( \omega \wedge \theta \right) \left( z_{1},...,z_{q+r}\right)
=\Gamma \left( Th\circ \rho ,h\circ \eta \right) z\left( \left( \omega
\wedge \theta \right) \left( z_{1},...,z_{q+r}\right) \right) \\ 
-\overset{q+r}{\underset{i=1}{\tsum }}\left( \omega \wedge \theta \right)
\left( \left( z_{1},...,\left[ z,z_{i}\right] _{F,h},...,z_{q+r}\right)
\right) \\ 
=\Gamma \left( Th\circ \rho ,h\circ \eta \right) z\left( \underset{\sigma
\left( q+1\right) <...<\sigma \left( q+r\right) }{\underset{\sigma \left(
1\right) <...<\sigma \left( q\right) }{\tsum }}sgn\left( \sigma \right)
\cdot \omega \left( z_{\sigma \left( 1\right) },...,z_{\sigma \left(
q\right) }\right) \right. \\ 
\qquad \left. \cdot \theta \left( z_{\sigma \left( q+1\right)
},...,z_{\sigma \left( q+r\right) }\right) \right) -\overset{q+r}{\underset{%
i=1}{\tsum }}\left( \omega \wedge \theta \right) \left( \left( z_{1},...,%
\left[ z,z_{i}\right] _{F,h},...,z_{q+r}\right) \right) \\ 
=\underset{\sigma \left( q+1\right) <...<\sigma \left( q+r\right) }{\underset%
{\sigma \left( 1\right) <...<\sigma \left( q\right) }{\tsum }}sgn\left(
\sigma \right) \cdot \Gamma \left( Th\circ \rho ,h\circ \eta \right) z\left(
\omega \left( z_{\sigma \left( 1\right) },...,z_{\sigma \left( q\right)
}\right) \right) \\ 
\qquad \cdot \theta \left( z_{\sigma \left( q+1\right) },...,z_{\sigma
\left( q+r\right) }\right) +\underset{\sigma \left( q+1\right) <...<\sigma
\left( q+r\right) }{\underset{\sigma \left( 1\right) <...<\sigma \left(
q\right) }{\tsum }}sgn\left( \sigma \right) \cdot \omega \left( z_{\sigma
\left( 1\right) },...,z_{\sigma \left( q\right) }\right) \\ 
\qquad \cdot \Gamma \left( Th\circ \rho ,h\circ \eta \right) z\left( \theta
\left( z_{\sigma \left( q+1\right) },...,z_{\sigma \left( q+r\right)
}\right) \right) -\underset{\sigma \left( q+1\right) <...<\sigma \left(
q+r\right) }{\underset{\sigma \left( 1\right) <...<\sigma \left( q\right) }{%
\tsum }}sgn\left( \sigma \right) \\ 
\qquad \cdot \overset{q}{\underset{i=1}{\tsum }}\omega \left( z_{\sigma
\left( 1\right) },...,\left[ z,z_{\sigma \left( i\right) }\right]
_{F,h},...,z_{\sigma \left( q\right) }\right) \cdot \theta \left( z_{\sigma
\left( q+1\right) },...,z_{\sigma \left( q+r\right) }\right)%
\end{array}%
\end{equation*}%
\begin{equation*}
\begin{array}{l}
-\underset{\sigma \left( q+1\right) <...<\sigma \left( q+r\right) }{\underset%
{\sigma \left( 1\right) <...<\sigma \left( q\right) }{\tsum }}sgn\left(
\sigma \right) \overset{q+r}{\underset{i=q+1}{\tsum }}\omega \left(
z_{\sigma \left( 1\right) },...,z_{\sigma \left( q\right) }\right) \\ 
\qquad \cdot \theta \left( z_{\sigma \left( q+1\right) },...,\left[
z,z_{\sigma \left( i\right) }\right] _{F,h},...,z_{\sigma \left( q+r\right)
}\right) \\ 
=\underset{\sigma \left( q+1\right) <...<\sigma \left( q+r\right) }{\underset%
{\sigma \left( 1\right) <...<\sigma \left( q\right) }{\tsum }}sgn\left(
\sigma \right) L_{z}\omega \left( z_{\sigma \left( 1\right) },...,\left[
z,z_{\sigma \left( i\right) }\right] _{F,h},...,z_{\sigma \left( q\right)
}\right) \\ 
\qquad \cdot \theta \left( z_{\sigma \left( q+1\right) },...,z_{\sigma
\left( q+r\right) }\right) +\underset{\sigma \left( q+1\right) <...<\sigma
\left( q+r\right) }{\underset{\sigma \left( 1\right) <...<\sigma \left(
q\right) }{\tsum }}sgn\left( \sigma \right) \overset{q+r}{\underset{i=q+1}{%
\tsum }}\omega \left( z_{\sigma \left( 1\right) },...,z_{\sigma \left(
q\right) }\right) \\ 
\qquad \cdot L_{z}\theta \left( z_{\sigma \left( q+1\right) },...,\left[
z,z_{\sigma \left( i\right) }\right] _{F,h},...,z_{\sigma \left( q+r\right)
}\right) \\ 
=\left( L_{z}\omega \wedge \theta +\omega \wedge L_{z}\theta \right) \left(
z_{1},...,z_{q+r}\right)%
\end{array}%
\end{equation*}%
it results the conclusion of the theorem. \hfill \emph{q.e.d.}

\textbf{Definition 2.6 }For any $z\in \Gamma \left( F,\nu ,N\right) $, the $%
\mathcal{F}\left( N\right) $-multilinear application 
\begin{equation*}
\begin{array}{rcl}
\Lambda \left( F,\nu ,N\right) & ^{\underrightarrow{\ \ i_{z}\ \ }} & 
\Lambda \left( F,\nu ,N\right) \\ 
\Lambda ^{q}\left( F,\nu ,N\right) \ni \omega & \longmapsto & i_{z}\omega
\in \Lambda ^{q-1}\left( F,\nu ,N\right) ,%
\end{array}%
\end{equation*}%
where 
\begin{equation*}
\begin{array}{c}
i_{z}\omega \left( z_{2},...,z_{q}\right) =\omega \left(
z,z_{2},...,z_{q}\right) ,%
\end{array}%
\end{equation*}%
for any $z_{2},...,z_{q}\in \Gamma \left( F,\nu ,N\right) $, will be called
the \emph{interior product associated to the section}~$z.$\bigskip

For any $f\in \mathcal{F}\left( N\right) $, we define $\ i_{z}f=0.$

\textit{Remark 2.4}\textbf{\ }If $z\in \Gamma \left( F,\nu ,N\right) ,$ $%
\omega \in $\ $\Lambda ^{p}\left( F,\nu ,N\right) $ and $U$\ is an open
subset of $N$\ such that $z_{|U}=0$\ or $\omega _{|U}=0,$\ then $\left(
i_{z}\omega \right) _{|U}=0.$

\textbf{Theorem 2.5 }\emph{If }$z\in \Gamma \left( F,\nu ,N\right) $\emph{,
then for any }$\omega \in $\emph{\ }$\Lambda ^{q}\left( F,\nu ,N\right) $%
\emph{\ and }$\theta \in $\emph{\ }$\Lambda ^{r}\left( F,\nu ,N\right) $%
\emph{\ we obtain} 
\begin{equation*}
\begin{array}{c}
i_{z}\left( \omega \wedge \theta \right) =i_{z}\omega \wedge \theta +\left(
-1\right) ^{q}\omega \wedge i_{z}\theta .%
\end{array}%
\leqno(2.10)
\end{equation*}

\emph{Proof. }Let $z_{1},...,z_{q+r}\in \Gamma \left( F,\nu ,N\right) $ be
arbitrary. We observe that 
\begin{equation*}
\begin{array}{l}
i_{z_{1}}\left( \omega \wedge \theta \right) \left( z_{2},...,z_{q+r}\right)
=\left( \omega \wedge \theta \right) \left( z_{1},z_{2},...,z_{q+r}\right)
\\ 
=\underset{\sigma \left( q+1\right) <...<\sigma \left( q+r\right) }{\underset%
{\sigma \left( 1\right) <...<\sigma \left( q\right) }{\tsum }}sgn\left(
\sigma \right) \cdot \omega \left( z_{\sigma \left( 1\right) },...,z_{\sigma
\left( q\right) }\right) \cdot \theta \left( z_{\sigma \left( q+1\right)
},...,z_{\sigma \left( q+r\right) }\right) \\ 
=\underset{\sigma \left( q+1\right) <...<\sigma \left( q+r\right) }{\underset%
{1=\sigma \left( 1\right) <\sigma \left( 2\right) <...<\sigma \left(
q\right) }{\tsum }}sgn\left( \sigma \right) \cdot \omega \left(
z_{1},z_{\sigma \left( 2\right) },...,z_{\sigma \left( q\right) }\right)
\cdot \theta \left( z_{\sigma \left( q+1\right) },...,z_{\sigma \left(
q+r\right) }\right) \\ 
+\underset{1=\sigma \left( q+1\right) <\sigma \left( q+2\right) <...<\sigma
\left( q+r\right) }{\underset{\sigma \left( 1\right) <...<\sigma \left(
q\right) }{\tsum }}sgn\left( \sigma \right) \cdot \omega \left( z_{\sigma
\left( 1\right) },...,z_{\sigma \left( q\right) }\right) \cdot \theta \left(
z_{1},z_{\sigma \left( q+2\right) },...,z_{\sigma \left( q+r\right) }\right)
\\ 
=\underset{\sigma \left( q+1\right) <...<\sigma \left( q+r\right) }{\underset%
{\sigma \left( 2\right) <...<\sigma \left( q\right) }{\tsum }}sgn\left(
\sigma \right) \cdot i_{z_{1}}\omega \left( z_{\sigma \left( 2\right)
},...,z_{\sigma \left( q\right) }\right) \cdot \theta \left( z_{\sigma
\left( q+1\right) },...,z_{\sigma \left( q+r\right) }\right) \\ 
+\underset{\sigma \left( q+2\right) <...<\sigma \left( q+r\right) }{\underset%
{\sigma \left( 1\right) <...<\sigma \left( q\right) }{\tsum }}sgn\left(
\sigma \right) \cdot \omega \left( z_{\sigma \left( 1\right) },...,z_{\sigma
\left( q\right) }\right) \cdot i_{z_{1}}\theta \left( z_{\sigma \left(
q+2\right) },...,z_{\sigma \left( q+r\right) }\right) .%
\end{array}%
\end{equation*}

In the second sum, we have the permutation%
\begin{equation*}
\sigma =\left( 
\begin{array}{ccccccc}
1 & ... & q & q+1 & q+2 & ... & q+r \\ 
\sigma \left( 1\right) & ... & \sigma \left( q\right) & 1 & \sigma \left(
q+2\right) & ... & \sigma \left( q+r\right)%
\end{array}%
\right) .
\end{equation*}

We observe that $\sigma =\tau \circ \tau ^{\prime }$, where%
\begin{equation*}
\tau =\left( 
\begin{array}{ccccccc}
1 & 2 & ... & q+1 & q+2 & ... & q+r \\ 
1 & \sigma \left( 1\right) & ... & \sigma \left( q\right) & \sigma \left(
q+2\right) & ... & \sigma \left( q+r\right)%
\end{array}%
\right)
\end{equation*}%
and 
\begin{equation*}
\tau ^{\prime }=\left( 
\begin{array}{cccccccc}
1 & 2 & ... & q & q+1 & q+2 & ... & q+r \\ 
2 & 3 & ... & q+1 & 1 & q+2 & ... & q+r%
\end{array}%
\right) .
\end{equation*}

Since $\tau \left( 2\right) <...<\tau \left( q+1\right) $ and $\tau ^{\prime
}$ has $q$ inversions, it results that 
\begin{equation*}
sgn\left( \sigma \right) =\left( -1\right) ^{q}\cdot sgn\left( \tau \right) .
\end{equation*}

Therefore, 
\begin{equation*}
\begin{array}{l}
i_{z_{1}}\left( \omega \wedge \theta \right) \left( z_{2},...,z_{q+r}\right)
=\left( i_{z_{1}}\omega \wedge \theta \right) \left( z_{2},...,z_{q+r}\right)
\\ 
+\left( -1\right) ^{q}\underset{\tau \left( q+2\right) <...<\tau \left(
q+r\right) }{\underset{\tau \left( 2\right) <...<\tau \left( q\right) }{%
\tsum }}sgn\left( \tau \right) \cdot \omega \left( z_{\tau \left( 2\right)
},...,z_{\tau \left( q\right) }\right) \cdot i_{z_{1}}\theta \left( z_{\tau
\left( q+2\right) },...,z_{\tau \left( q+r\right) }\right) \\ 
=\left( i_{z_{1}}\omega \wedge \theta \right) \left(
z_{2},...,z_{q+r}\right) +\left( -1\right) ^{q}\left( \omega \wedge
i_{z_{1}}\theta \right) \left( z_{2},...,z_{q+r}\right) .%
\end{array}%
\end{equation*}

\hfill \emph{q.e.d.}

\textbf{Theorem 2.6 }\emph{For any }$z,v\in \Gamma \left( F,\nu ,N\right) $%
\emph{\ we obtain}%
\begin{equation*}
\begin{array}{c}
L_{v}\circ i_{z}-i_{z}\circ L_{v}=i_{\left[ z,v\right] _{F,h}}.%
\end{array}%
\leqno(2.11)
\end{equation*}

\emph{Proof.} Let $\omega \in $\emph{\ }$\Lambda ^{q}\left( F,\nu ,N\right) $
be arbitrary. Since 
\begin{equation*}
\begin{array}{l}
i_{z}\left( L_{v}\omega \right) \left( z_{2},...z_{q}\right) =L_{v}\omega
\left( z,z_{2},...z_{q}\right) \\ 
=\Gamma \left( Th\circ \rho ,h\circ \eta \right) v\left( \omega \left(
z,z_{2},...,z_{q}\right) \right) -\omega \left( \left[ v,z\right]
_{F,h},z_{2},...,z_{q}\right) \\ 
-\overset{q}{\underset{i=2}{\tsum }}\omega \left( \left( z,z_{2},...,\left[
v,z_{i}\right] _{F,h},...,z_{q}\right) \right) \\ 
=\Gamma \left( Th\circ \rho ,h\circ \eta \right) v\left( i_{z}\omega \left(
z_{2},...,z_{q}\right) \right) -\overset{q}{\underset{i=2}{\tsum }}%
i_{z}\omega \left( z_{2},...,\left[ v,z_{i}\right] _{F,h},...,z_{q}\right)
\\ 
-i_{\left[ v,z\right] _{F,h}}\left( z_{2},...,z_{q}\right) =\left(
L_{v}\left( i_{z}\omega \right) -i_{\left[ v,z\right] _{F,h}}\right) \left(
z_{2},...,z_{q}\right) ,%
\end{array}%
\end{equation*}%
for any $z_{2},...,z_{q}\in \Gamma \left( F,\nu ,N\right) $ it result the
conclusion of the theorem.\hfill \emph{q.e.d.}

\textbf{Definition 2.7 }If $f\in \mathcal{F}\left( N\right) $ and $z\in
\Gamma \left( F,\nu ,N\right) ,$ then \emph{the exterior differential
operator} is defined by \textbf{\ }%
\begin{equation*}
\begin{array}{c}
d^{F}f\left( z\right) =\Gamma \left( Th\circ \rho ,h\circ \eta \right)
\left( z\right) f.%
\end{array}%
\leqno(2.12)
\end{equation*}

In the particular case of Lie algebroids, $\left( \eta ,h\right) =\left(
Id_{M},Id_{M}\right) ,$ then \emph{the exterior differential operator} is
defined by 
\begin{equation*}
\begin{array}{c}
d^{F}f\left( z\right) =\Gamma \left( \rho ,Id_{M}\right) \left( z\right) f.%
\end{array}%
\leqno(2.12^{\prime })
\end{equation*}

In addition, if $\rho =Id_{TM},$ then \emph{the exterior differential
operator} is defined by 
\begin{equation*}
\begin{array}{c}
d^{TM}f\left( z\right) =z\left( f\right) .%
\end{array}%
\leqno(2.12^{\prime \prime })
\end{equation*}

\textbf{Theorem 2.7 }\emph{The }$\mathcal{F}\left( N\right) $\emph{%
-multilinear application }%
\begin{equation*}
\begin{array}{c}
\begin{array}{ccc}
\Lambda ^{q}\mathbf{\ }\left( F,\nu ,N\right) & ^{\underrightarrow{\,\
d^{F}\,\ }} & \Lambda ^{q+1}\mathbf{\ }\left( F,\nu ,N\right) \\ 
\omega & \longmapsto & d\omega%
\end{array}%
\end{array}%
\end{equation*}%
\emph{defined by} 
\begin{equation*}
\begin{array}{l}
d^{F}\omega \left( z_{0},z_{1},...,z_{q}\right) =\overset{q}{\underset{i=0}{%
\tsum }}\left( -1\right) ^{i}\Gamma \left( Th\circ \rho ,h\circ \eta \right)
z_{i}\left( \omega \left( z_{0},z_{1},...,\hat{z}_{i},...,z_{q}\right)
\right) \\ 
~\ \ \ \ \ \ \ \ \ \ \ \ \ \ \ \ \ \ \ \ \ \ \ \ \ \ +\underset{i<j}{\tsum }%
\left( -1\right) ^{i+j}\omega \left( \left[ z_{i},z_{j}\right]
_{F,h},z_{0},z_{1},...,\hat{z}_{i},...,\hat{z}_{j},...,z_{q}\right) ,%
\end{array}%
\leqno(2.13)
\end{equation*}%
\emph{for any }$z_{0},z_{1},...,z_{q}\in \Gamma \left( F,\nu ,N\right) ,$ 
\emph{is unique with the following property:}%
\begin{equation*}
\begin{array}{c}
L_{z}=d^{F}\circ i_{z}+i_{z}\circ d^{F},~\forall z\in \Gamma \left( F,\nu
,N\right) .%
\end{array}%
\leqno(2.14)
\end{equation*}

This $\mathcal{F}\left( N\right) $-multilinear application will be called%
\emph{\ the exterior differentiation ope\-ra\-tor for the exterior
differential algebra of the generalized Lie algebroid }$((F,\nu
,N),[,]_{F,h},(\rho ,\eta )).$

In the particular case of Lie algebroids, $\left( \eta ,h\right) =\left(
Id_{M},Id_{M}\right) ,$ then \emph{the exterior differentiation ope\-ra\-tor
for the exterior differential algebra of the Lie algebroid }$((F,\nu
,M),[,]_{F},(\rho ,Id_{M}))$ is defined by 
\begin{equation*}
\begin{array}{l}
d^{F}\omega \left( z_{0},z_{1},...,z_{q}\right) =\overset{q}{\underset{i=0}{%
\tsum }}\left( -1\right) ^{i}\Gamma \left( \rho ,Id_{M}\right) z_{i}\left(
\omega \left( z_{0},z_{1},...,\hat{z}_{i},...,z_{q}\right) \right) \\ 
~\ \ \ \ \ \ \ \ \ \ \ \ \ \ \ \ \ \ \ \ \ \ \ \ \ \ +\underset{i<j}{\tsum }%
\left( -1\right) ^{i+j}\omega \left( \left[ z_{i},z_{j}\right]
_{F},z_{0},z_{1},...,\hat{z}_{i},...,\hat{z}_{j},...,z_{q}\right) ,%
\end{array}%
\leqno(2.13^{\prime })
\end{equation*}

In addition, if $\rho =Id_{TM},$ then \emph{the exterior differentiation
ope\-ra\-tor for the exterior differential algebra of the standard Lie
algebroid }$((TM,\nu ,M),[,]_{TM},(Id_{TM},Id_{M}))$ is defined by 
\begin{equation*}
\begin{array}{l}
d^{TM}\omega \left( z_{0},z_{1},...,z_{q}\right) =\overset{q}{\underset{i=0}{%
\tsum }}\left( -1\right) ^{i}z_{i}\left( \omega \left( z_{0},z_{1},...,\hat{z%
}_{i},...,z_{q}\right) \right) \\ 
~\ \ \ \ \ \ \ \ \ \ \ \ \ \ \ \ \ \ \ \ \ \ \ \ \ \ +\underset{i<j}{\tsum }%
\left( -1\right) ^{i+j}\omega \left( \left[ z_{i},z_{j}\right]
_{TM},z_{0},z_{1},...,\hat{z}_{i},...,\hat{z}_{j},...,z_{q}\right) ,%
\end{array}%
\leqno(2.13^{\prime \prime })
\end{equation*}

\emph{Proof. }We verify the property $\left( 2.13\right) $ Since 
\begin{equation*}
\begin{array}{l}
\left( i_{z_{0}}\circ d^{F}\right) \omega \left( z_{1},...,z_{q}\right)
=d\omega \left( z_{0},z_{1},...,z_{q}\right) \\ 
=\overset{q}{\underset{i=0}{\tsum }}\left( -1\right) ^{i}\Gamma \left(
Th\circ \rho ,h\circ \eta \right) z_{i}\left( \omega \left( z_{0},z_{1},...,%
\hat{z}_{i},...,z_{q}\right) \right) \\ 
+\underset{0\leq i<j}{\tsum }\left( -1\right) ^{i+j}\omega \left( \left[
z_{i},z_{j}\right] _{F,h},z_{0},z_{1},...,\hat{z}_{i},...,\hat{z}%
_{j},...,z_{q}\right) \\ 
=\Gamma \left( Th\circ \rho ,h\circ \eta \right) z_{0}\left( \omega \left(
z_{1},...,z_{q}\right) \right) \\ 
+\overset{q}{\underset{i=1}{\tsum }}\left( -1\right) ^{i}\Gamma \left(
Th\circ \rho ,h\circ \eta \right) z_{i}\left( \omega \left( z_{0},z_{1},...,%
\hat{z}_{i},...,z_{q}\right) \right) \\ 
+\underset{i=1}{\overset{q}{\tsum }}\left( -1\right) ^{i}\omega \left( \left[
z_{0},z_{i}\right] _{F,h},z_{1},...,\hat{z}_{i},...,z_{q}\right) \\ 
+\underset{1\leq i<j}{\tsum }\left( -1\right) ^{i+j}\omega \left( \left[
z_{i},z_{j}\right] _{F,h},z_{0},z_{1},...,\hat{z}_{i},...,\hat{z}%
_{j},...,z_{q}\right) \\ 
=\Gamma \left( Th\circ \rho ,h\circ \eta \right) z_{0}\left( \omega \left(
z_{1},...,z_{q}\right) \right) \\ 
-\underset{i=1}{\overset{q}{\tsum }}\omega \left( z_{1},...,\left[
z_{0},z_{i}\right] _{F,h},...,z_{q}\right) \\ 
-\overset{q}{\underset{i=1}{\tsum }}\left( -1\right) ^{i-1}\Gamma \left(
Th\circ \rho ,h\circ \eta \right) z_{i}\left( i_{z_{0}}\omega \left( \left(
z_{1},...,\hat{z}_{i},...,z_{q}\right) \right) \right) \\ 
-\underset{1\leq i<j}{\tsum }\left( -1\right) ^{i+j-2}i_{z_{0}}\omega \left(
\left( \left[ z_{i},z_{j}\right] _{F,h},z_{1},...,\hat{z}_{i},...,\hat{z}%
_{j},...,z_{q}\right) \right) \\ 
=\left( L_{z_{0}}-d^{F}\circ i_{z_{0}}\right) \omega \left(
z_{1},...,z_{q}\right) ,%
\end{array}%
\end{equation*}%
for any $z_{0},z_{1},...,z_{q}\in \Gamma \left( F,\nu ,N\right) $ it results
that the property $\left( 2.13\right) $ is satisfied.

In the following, we verify the uniqueness of the operator $d^{F}.$

Let $d^{\prime F}$ be an another exterior differentiation operator
satisfying the property $\left( 2.13\right) .$

Let $S=\left\{ q\in \mathbb{N}:d^{F}\omega =d^{\prime F}\omega ,~\forall
\omega \in \Lambda ^{q}\left( F,\nu ,N\right) \right\} $ be.

Let $z\in \Gamma \left( F,\nu ,N\right) $ be arbitrary.

We observe that $\left( 2.13\right) $ is equivalent with 
\begin{equation*}
\begin{array}{c}
i_{z}\circ \left( d^{F}-d^{\prime F}\right) +\left( d^{F}-d^{\prime
F}\right) \circ i_{z}=0.%
\end{array}%
\leqno(1)
\end{equation*}

Since $i_{z}f=0,$ for any $f\in \mathcal{F}\left( N\right) ,$ it results
that 
\begin{equation*}
\begin{array}{c}
\left( \left( d^{F}-d^{\prime F}\right) f\right) \left( z\right) =0,~\forall
f\in \mathcal{F}\left( N\right) .%
\end{array}%
\end{equation*}

Therefore, we obtain that 
\begin{equation*}
\begin{array}{c}
0\in S.%
\end{array}%
\leqno(2)
\end{equation*}

In the following, we prove that 
\begin{equation*}
\begin{array}{c}
q\in S\Longrightarrow q+1\in S%
\end{array}%
\leqno(3)
\end{equation*}

Let $\omega \in \Lambda ^{p+1}\left( F,\nu ,N\right) $ be arbitrary$.$ Since 
$i_{z}\omega \in \Lambda ^{q}\left( F,\nu ,N\right) $, using the equality $%
\left( 1\right) $, it results that 
\begin{equation*}
i_{z}\circ \left( d^{F}-d^{\prime F}\right) \omega =0.
\end{equation*}

We obtain that, $\left( \left( d^{F}-d^{\prime F}\right) \omega \right)
\left( z_{0},z_{1},...,z_{q}\right) =0,$ for any $z_{1},...,z_{q}\in \Gamma
\left( F,\nu ,N\right) .$

Therefore $d^{F}\omega =d^{\prime F}\omega ,$ namely $q+1\in S.$

Using the \textit{Peano's Axiom }and the affirmations $\left( 2\right) $ and 
$\left( 3\right) $ it results that $S=\mathbb{N}.$

Therefore, the uniqueness is verified.\hfill \emph{q.e.d.}\medskip

Note that if $\omega =\omega _{\alpha _{1}...\alpha _{q}}t^{\alpha
_{1}}\wedge ...\wedge t^{\alpha _{q}}\in \Lambda ^{q}\left( F,\nu ,N\right) $%
, then 
\begin{eqnarray*}
d^{F}\omega \left( t_{\alpha _{0}},t_{\alpha _{1}},...,t_{\alpha
_{q}}\right) &=&\overset{q}{\underset{i=0}{\tsum }}\left( -1\right)
^{i}\theta _{\alpha _{i}}^{\tilde{k}}\frac{\partial \omega _{\alpha _{0},...,%
\widehat{\alpha _{i}}...\alpha _{q}}}{\partial \varkappa ^{\tilde{k}}} \\%
[1mm]
&&+\underset{i<j}{\tsum }\left( -1\right) ^{i+j}L_{\alpha _{i}\alpha
_{j}}^{\alpha }\cdot \omega _{\alpha ,\alpha _{0},...,\widehat{\alpha _{i}}%
,...,\widehat{\alpha _{j}},...,\alpha _{q}}.
\end{eqnarray*}

Therefore, we obtain%
\begin{equation*}
\begin{array}{l}
d^{F}\omega =\left( \overset{q}{\underset{i=0}{\tsum }}\left( -1\right)
^{i}\theta _{\alpha _{i}}^{\tilde{k}}\displaystyle\frac{\partial \omega
_{\alpha _{0},...,\widehat{\alpha _{i}}...\alpha _{q}}}{\partial \varkappa ^{%
\tilde{k}}}\right. \\ 
\left. +\underset{i<j}{\tsum }\left( -1\right) ^{i+j}L_{\alpha _{i}\alpha
_{j}}^{\alpha }\cdot \omega _{\alpha ,\alpha _{0},...,\widehat{\alpha _{i}}%
,...,\widehat{\alpha _{j}},...,\alpha _{q}}\right) t^{\alpha _{0}}\wedge
t^{\alpha _{1}}\wedge ...\wedge t^{\alpha _{q}}.%
\end{array}%
\leqno(2.15)
\end{equation*}

\emph{Remark 2.5 }If $d^{F}$ is the exterior differentiation operator for
the generalized Lie algebroid 
\begin{equation*}
\left( \left( F,\nu ,N\right) ,\left[ ,\right] _{F,h},\left( \rho ,\eta
\right) \right) ,
\end{equation*}%
$\omega \in $\emph{\ }$\Lambda ^{q}\left( F,\nu ,N\right) $\ and $U$\ is an
open subset of $N$\ such that $\omega _{|U}=0,$\ then $\left( d^{F}\omega
\right) _{|U}=0.$

\textbf{Theorem 2.8} \emph{The exterior differentiation operator }$d^{F}$%
\emph{\ given by the previous theorem has the following properties:}\medskip

\noindent 1. \emph{\ For any }$\omega \in $\emph{\ }$\Lambda ^{q}\left(
F,\nu ,N\right) $\emph{\ and }$\theta \in $\emph{\ }$\Lambda ^{r}\left(
F,\nu ,N\right) $\emph{\ we obtain }%
\begin{equation*}
\begin{array}{c}
d^{F}\left( \omega \wedge \theta \right) =d^{F}\omega \wedge \theta +\left(
-1\right) ^{q}\omega \wedge d^{F}\theta .%
\end{array}%
\leqno(2.16)
\end{equation*}

\noindent 2.\emph{\ For any }$z\in \Gamma \left( F,\nu ,N\right) $ we obtain 
\begin{equation*}
\begin{array}{c}
L_{z}\circ d^{F}=d^{F}\circ L_{z}.%
\end{array}%
\leqno(2.17)
\end{equation*}

\noindent3.\emph{\ }$d^{F}\circ d^{F}=0.$

\emph{Proof.}\newline
1. Let $S=\left\{ q\in \mathbb{N}:d^{F}\left( \omega \wedge \theta \right)
=d^{F}\omega \wedge \theta +\left( -1\right) ^{q}\omega \wedge d^{F}\theta
,~\forall \omega \in \Lambda ^{q}\left( F,\nu ,N\right) \right\} $ be. Since 
\begin{equation*}
\begin{array}{l}
d^{F}\left( f\wedge \theta \right) \left( z,v\right) =d^{F}\left( f\cdot
\theta \right) \left( z,v\right) \\ 
=\Gamma \left( Th\circ \rho ,h\circ \eta \right) z\left( f\omega \left(
v\right) \right) -\Gamma \left( Th\circ \rho ,h\circ \eta \right) v\left(
f\omega \left( z\right) \right) -f\omega \left( \left[ z,v\right]
_{F,h}\right) \\ 
=\Gamma \left( Th\circ \rho ,h\circ \eta \right) z\left( f\right) \cdot
\omega \left( v\right) +f\cdot \Gamma \left( Th\circ \rho ,h\circ \eta
\right) z\left( \omega \left( v\right) \right) \\ 
-\Gamma \left( Th\circ \rho ,h\circ \eta \right) v\left( f\right) \cdot
\omega \left( z\right) -f\cdot \Gamma \left( Th\circ \rho ,h\circ \eta
\right) v\left( \omega \left( z\right) \right) -f\omega \left( \left[ z,v%
\right] _{F,h}\right) \\ 
=d^{F}f\left( z\right) \cdot \omega \left( v\right) -d^{F}f\left( v\right)
\cdot \omega \left( z\right) +f\cdot d^{F}\omega \left( z,v\right) \\ 
=\left( d^{F}f\wedge \omega \right) \left( z,v\right) +\left( -1\right)
^{0}f\cdot d^{F}\omega \left( z,v\right) \\ 
=\left( d^{F}f\wedge \omega \right) \left( z,v\right) +\left( -1\right)
^{0}\left( f\wedge d^{F}\omega \right) \left( z,v\right) ,~\forall z,v\in
\Gamma \left( F,\nu ,N\right) ,%
\end{array}%
\end{equation*}%
it results that 
\begin{equation*}
\begin{array}{c}
0\in S.%
\end{array}%
\leqno(1.1)
\end{equation*}

In the following we prove that 
\begin{equation*}
\begin{array}{c}
q\in S\Longrightarrow q+1\in S.%
\end{array}%
\leqno(1.2)
\end{equation*}

Without restricting the generality, we consider that $\theta \in \Lambda
^{r}\left( F,\nu ,N\right) .$ Since 
\begin{equation*}
\begin{array}{l}
d^{F}\left( \omega \wedge \theta \right) \left(
z_{0},z_{1},...,z_{q+r}\right) =i_{z_{0}}\circ d^{F}\left( \omega \wedge
\theta \right) \left( z_{1},...,z_{q+r}\right) \\ 
=L_{z_{0}}\left( \omega \wedge \theta \right) \left(
z_{1},...,z_{q+r}\right) -d^{F}\circ i_{z_{0}}\left( \omega \wedge \theta
\right) \left( z_{1},...,z_{q+r}\right) \\ 
=\left( L_{z_{0}}\omega \wedge \theta +\omega \wedge L_{z_{0}}\theta \right)
\left( z_{1},...,z_{q+r}\right) \\ 
-\left[ d^{F}\circ \left( i_{z_{0}}\omega \wedge \theta +\left( -1\right)
^{q}\omega \wedge i_{z_{0}}\theta \right) \right] \left(
z_{1},...,z_{q+r}\right) \\ 
=\left( L_{z_{0}}\omega \wedge \theta +\omega \wedge L_{z_{0}}\theta -\left(
d^{F}\circ i_{z_{0}}\omega \right) \wedge \theta \right) \left(
z_{1},...,z_{q+r}\right) \\ 
-\left( \left( -1\right) ^{q-1}i_{z_{0}}\omega \wedge d^{F}\theta +\left(
-1\right) ^{q}d^{F}\omega \wedge i_{z_{0}}\theta \right) \left(
z_{1},...,z_{q+r}\right) \\ 
-\left( -1\right) ^{2q}\omega \wedge d^{F}\circ i_{z_{0}}\theta \left(
z_{1},...,z_{q+r}\right) \\ 
=\left( \left( L_{z_{0}}\omega -d^{F}\circ i_{z_{0}}\omega \right) \wedge
\theta \right) \left( z_{1},...,z_{q+r}\right) \\ 
+\omega \wedge \left( L_{z_{0}}\theta -d^{F}\circ i_{z_{0}}\theta \right)
\left( z_{1},...,z_{q+r}\right) \\ 
+\left( \left( -1\right) ^{q}i_{z_{0}}\omega \wedge d^{F}\theta -\left(
-1\right) ^{q}d^{F}\omega \wedge i_{z_{0}}\theta \right) \left(
z_{1},...,z_{q+r}\right) \\ 
=\left[ \left( \left( i_{z_{0}}\circ d^{F}\right) \omega \right) \wedge
\theta +\left( -1\right) ^{q+1}d^{F}\omega \wedge i_{z_{0}}\theta \right]
\left( z_{1},...,z_{q+r}\right) \\ 
+\left[ \omega \wedge \left( \left( i_{z_{0}}\circ d^{F}\right) \theta
\right) +\left( -1\right) ^{q}i_{z_{0}}\omega \wedge d^{F}\theta \right]
\left( z_{1},...,z_{q+r}\right) \\ 
=\left[ i_{z_{0}}\left( d^{F}\omega \wedge \theta \right) +\left( -1\right)
^{q}i_{z_{0}}\left( \omega \wedge d^{F}\theta \right) \right] \left(
z_{1},...,z_{q+r}\right) \\ 
=\left[ d^{F}\omega \wedge \theta +\left( -1\right) ^{q}\omega \wedge
d^{F}\theta \right] \left( z_{1},...,z_{q+r}\right) ,%
\end{array}%
\end{equation*}%
for any $z_{0},z_{1},...,z_{q+r}\in \Gamma \left( F,\nu ,N\right) $, it
results $\left( 1.2\right) .$

Using the \textit{Peano's Axiom }and the affirmations $\left( 1.1\right) $
and $\left( 1.2\right) $ it results that $S=\mathbb{N}.$

Therefore, it results the conclusion of affirmation 1.\medskip

2. Let $z\in \Gamma \left( F,\nu ,N\right) $ be arbitrary.

Let $S=\left\{ q\in \mathbb{N}:\left( L_{z}\circ d^{F}\right) \omega =\left(
d^{F}\circ L_{z}\right) \omega ,~\forall \omega \in \Lambda ^{q}\left( F,\nu
,N\right) \right\} $ be.

Let $f\in \mathcal{F}\left( N\right) $ be arbitrary. Since 
\begin{equation*}
\begin{array}{l}
\left( d^{F}\circ L_{z}\right) f\left( v\right) =i_{v}\circ \left(
d^{F}\circ L_{z}\right) f=\left( i_{v}\circ d^{F}\right) \circ L_{z}f \\ 
=\left( L_{v}\circ L_{z}\right) f-\left( \left( d^{F}\circ i_{v}\right)
\circ L_{z}\right) f \\ 
=\left( L_{v}\circ L_{z}\right) f-L_{\left[ z,v\right] _{F,h}}f+d^{F}\circ
i_{\left[ z,v\right] _{F,h}}f-d^{F}\circ L_{z}\left( i_{v}f\right) \\ 
=\left( L_{v}\circ L_{z}\right) f-L_{\left[ z,v\right] _{F,h}}f+d^{F}\circ
i_{\left[ z,v\right] _{F,h}}f-0 \\ 
=\left( L_{v}\circ L_{z}\right) f-L_{\left[ z,v\right] _{F,h}}f+d^{F}\circ
i_{\left[ z,v\right] _{F,h}}f-L_{z}\circ d^{F}\left( i_{v}f\right) \\ 
=\left( L_{z}\circ i_{v}\right) \left( d^{F}f\right) -L_{\left[ z,v\right]
_{F,h}}f+d^{F}\circ i_{\left[ z,v\right] _{F,h}}f \\ 
=\left( i_{v}\circ L_{z}\right) \left( d^{F}f\right) +L_{\left[ z,v\right]
_{F,h}}f-L_{\left[ z,v\right] _{F,h}}f \\ 
=i_{v}\circ \left( L_{z}\circ d^{F}\right) f=\left( L_{z}\circ d^{F}\right)
f\left( v\right) ,~\forall v\in \Gamma \left( F,\nu ,N\right) ,%
\end{array}%
\end{equation*}%
it results that 
\begin{equation*}
\begin{array}{c}
0\in S.%
\end{array}%
\leqno(2.1)
\end{equation*}

In the following we prove that 
\begin{equation*}
\begin{array}{c}
q\in S\Longrightarrow q+1\in S.%
\end{array}%
\leqno(2.2)
\end{equation*}

Let $\omega \in \Lambda ^{q}\left( F,\nu ,N\right) $ be arbitrary. Since 
\begin{equation*}
\begin{array}{l}
\left( d^{F}\circ L_{z}\right) \omega \left( z_{0},z_{1},...,z_{q}\right)
=i_{z_{0}}\circ \left( d^{F}\circ L_{z}\right) \omega \left(
z_{1},...,z_{q}\right) \vspace*{1mm} \\ 
=\left( i_{z_{0}}\circ d^{F}\right) \circ L_{z}\omega \left(
z_{1},...,z_{q}\right) \vspace*{1mm} \\ 
=\left[ \left( L_{z_{0}}\circ L_{z}\right) \omega -\left( \left( d^{F}\circ
i_{z_{0}}\right) \circ L_{z}\right) \omega \right] \left(
z_{1},...,z_{q}\right) \vspace*{1mm} \\ 
=\left[ \left( L_{z_{0}}\circ L_{z}\right) \omega -L_{\left[ z,z_{0}\right]
_{F,h}}\omega \right] \left( z_{1},...,z_{q}\right) \vspace*{1mm} \\ 
+\left[ d^{F}\circ i_{\left[ z,z_{0}\right] _{F,h}}\omega -d^{F}\circ
L_{z}\left( i_{z_{0}}\omega \right) \right] \left( z_{1},...,z_{q}\right) 
\vspace*{1mm} \\ 
\overset{ip.}{=}\left[ \left( L_{z_{0}}\circ L_{z}\right) \omega -L_{\left[
z,z_{0}\right] _{F,h}}\omega \right] \left( z_{1},...,z_{q}\right) \vspace*{%
1mm} \\ 
+\left[ d^{F}\circ i_{\left[ z,z_{0}\right] _{F,h}}\omega -L_{z}\circ
d^{F}\left( i_{z_{0}}\omega \right) \right] \left( z_{1},...,z_{q}\right) 
\vspace*{1mm} \\ 
=\left[ \left( L_{z}\circ i_{z_{0}}\right) \left( d^{F}\omega \right) -L_{%
\left[ z,z_{0}\right] _{F,h}}\omega +d^{F}\circ i_{\left[ z,z_{0}\right]
_{F,h}}\omega \right] \left( z_{1},...,z_{q}\right) \vspace*{1mm} \\ 
=\left[ \left( i_{z_{0}}\circ L_{z}\right) \left( d^{F}\omega \right) +L_{%
\left[ z,z_{0}\right] _{F,h}}\omega -L_{\left[ z,z_{0}\right] _{F,h}}\omega %
\right] \left( z_{1},...,z_{q}\right) \vspace*{1mm} \\ 
=i_{z_{0}}\circ \left( L_{z}\circ d^{F}\right) \omega \left(
z_{1},...,z_{q}\right) \vspace*{1mm} \\ 
=\left( L_{z}\circ d^{F}\right) \omega \left( z_{0},z_{1},...,z_{q}\right)
,~\forall z_{0},z_{1},...,z_{q}\in \Gamma \left( F,\nu ,N\right) ,%
\end{array}%
\end{equation*}%
it results $\left( 2.2\right) .$

Using the \textit{Peano's Axiom }and the affirmations $\left( 2.1\right) $
and $\left( 2.2\right) $ it results that $S=\mathbb{N}.$

Therefore, it results the conclusion of affirmation 2.\medskip

3. It is remarked that 
\begin{equation*}
\begin{array}{l}
i_{z}\circ \left( d^{F}\circ d^{F}\right) =\left( i_{z}\circ d^{F}\right)
\circ d^{F}=L_{z}\circ d^{F}-\left( d^{F}\circ i_{z}\right) \circ d^{F} 
\vspace*{1mm} \\ 
=L_{z}\circ d^{F}-d^{F}\circ L_{z}+d^{F}\circ \left( d^{F}\circ i_{z}\right)
=\left( d^{F}\circ d^{F}\right) \circ i_{z},%
\end{array}%
\end{equation*}%
for any $z\in \Gamma \left( F,\nu ,N\right) .$

Let $\omega \in \Lambda ^{q}\left( F,\nu ,N\right) $ be arbitrary. Since 
\begin{equation*}
\begin{array}{l}
\left( d^{F}\circ d^{F}\right) \omega \left( z_{1},...,z_{q+2}\right)
=i_{z_{q+2}}\circ ...\circ i_{z_{1}}\circ \left( d^{F}\circ d^{F}\right)
\omega =...\vspace*{1mm} \\ 
=i_{z_{q+2}}\circ \left( d^{F}\circ d^{F}\right) \circ i_{z_{q+1}}\left(
\omega \left( z_{1},...,z_{q}\right) \right) \vspace*{1mm} \\ 
=i_{z_{q+2}}\circ \left( d^{F}\circ d^{F}\right) \left( 0\right) =0,~\forall
z_{1},...,z_{q+2}\in \Gamma \left( F,\nu ,N\right) ,%
\end{array}%
\end{equation*}%
it results the conclusion of affirmation 3. \hfill \emph{q.e.d.}

\textbf{Theorem 2.9} \emph{If }$d^{F}$ \emph{is the exterior differentiation
operator for the exterior differential} $\mathcal{F}(N)$\emph{-algebra} $%
(\Lambda (F,\nu ,N),+,\cdot ,\wedge )$, \emph{then we obtain the structure
equations of Maurer-Cartan type }%
\begin{equation*}
\begin{array}{c}
d^{F}t^{\alpha }=-\displaystyle\frac{1}{2}L_{\beta \gamma }^{\alpha
}t^{\beta }\wedge t^{\gamma },~\alpha \in \overline{1,p}%
\end{array}%
\leqno(\mathcal{C}_{1})
\end{equation*}%
\emph{and\ }%
\begin{equation*}
\begin{array}{c}
d^{F}\varkappa ^{\tilde{\imath}}=\theta _{\alpha }^{\tilde{\imath}}t^{\alpha
},~\tilde{\imath}\in \overline{1,n},%
\end{array}%
\leqno(\mathcal{C}_{2})
\end{equation*}%
\emph{where }$\left\{ t^{\alpha },\alpha \in \overline{1,p}\right\} ~$\emph{%
is the coframe of the vector bundle }$\left( F,\nu ,N\right) .$\bigskip

This equations will be called \emph{the structure equations of Maurer-Cartan
type associa\-ted to the generalized Lie algebroid }$\left( \left( F,\nu
,N\right) ,\left[ ,\right] _{F,h},\left( \rho ,\eta \right) \right) .$

In the particular case of Lie algebroids, $\left( \eta ,h\right) =\left(
Id_{M},Id_{M}\right) ,$ then \emph{the structure equations of Maurer-Cartan
type }become%
\begin{equation*}
\begin{array}{c}
d^{F}t^{\alpha }=-\displaystyle\frac{1}{2}L_{\beta \gamma }^{\alpha
}t^{\beta }\wedge t^{\gamma },~\alpha \in \overline{1,p}%
\end{array}%
\leqno(\mathcal{C}_{1}^{\prime })
\end{equation*}%
\emph{and\ }%
\begin{equation*}
\begin{array}{c}
d^{F}x^{i}=\rho _{\alpha }^{i}t^{\alpha },~i\in \overline{1,m}.%
\end{array}%
\leqno(\mathcal{C}_{2}^{\prime })
\end{equation*}

In the particular case of standard Lie algebroid, $\rho =Id_{TM},$ then 
\emph{the structure equations of Maurer-Cartan type }become%
\begin{equation*}
\begin{array}{c}
d^{TM}dx^{i}=0,~i\in \overline{1,m}%
\end{array}%
\leqno(\mathcal{C}_{1}^{\prime \prime })
\end{equation*}%
\emph{and\ }%
\begin{equation*}
\begin{array}{c}
d^{TM}x^{i}=dx^{i},~i\in \overline{1,m}.%
\end{array}%
\leqno(\mathcal{C}_{2}^{\prime \prime })
\end{equation*}

\emph{Proof.} Let $\alpha \in \overline{1,p}$ be arbitrary. Since%
\begin{equation*}
\begin{array}{c}
d^{F}t^{\alpha }\left( t_{\beta },t_{\gamma }\right) =-L_{\beta \gamma
}^{\alpha },~\forall \beta ,\gamma \in \overline{1,p}%
\end{array}%
\end{equation*}%
it results that 
\begin{equation*}
\begin{array}{c}
d^{F}t^{\alpha }=-\underset{\beta <\gamma }{\tsum }L_{\beta \gamma }^{\alpha
}t^{\beta }\wedge t^{\gamma }.%
\end{array}%
\leqno(1)
\end{equation*}

Since $L_{\beta \gamma }^{\alpha }=-L_{\gamma \beta }^{\alpha }$ and $%
t^{\beta }\wedge t^{\gamma }=-t^{\gamma }\wedge t^{\beta }$, for nay $\beta
,\gamma \in \overline{1,p},$ it results that 
\begin{equation*}
\begin{array}{c}
\underset{\beta <\gamma }{\tsum }L_{\beta \gamma }^{\alpha }t^{\beta }\wedge
t^{\gamma }=\displaystyle\frac{1}{2}L_{\beta \gamma }^{\alpha }t^{\beta
}\wedge t^{\gamma }%
\end{array}%
\leqno(2)
\end{equation*}

Using the equalities $\left( 1\right) $ and $\left( 2\right) $ it results
the structure equation $(\mathcal{C}_{1}).$

Let $\tilde{\imath}\in \overline{1,n}$ be arbitrarily. Since 
\begin{equation*}
\begin{array}{c}
d^{F}\varkappa ^{\tilde{\imath}}\left( t_{\alpha }\right) =\theta _{\alpha
}^{\tilde{\imath}},~\forall \alpha \in \overline{1,p}%
\end{array}%
\end{equation*}%
it results the structure equation $(\mathcal{C}_{2}).$\hfill \emph{q.e.d.}%
\bigskip

\bigskip

Let $\left( \left( F^{\prime },\nu ^{\prime },N^{\prime }\right) ,\left[ ,%
\right] _{F^{\prime },h^{\prime }},\left( \rho ^{\prime },\eta ^{\prime
}\right) \right) $ be an another generalized Lie algebroid.

In the category $\mathbf{GLA,}$ we defined (see $\left[ 1\right] $) the set
of morphisms of 
\begin{equation*}
\left( \left( F,\nu ,N\right) ,\left[ ,\right] _{F,h},\left( \rho ,\eta
\right) \right)
\end{equation*}%
source and 
\begin{equation*}
\left( \left( F^{\prime },\nu ^{\prime },N^{\prime }\right) ,\left[ ,\right]
_{F^{\prime },h^{\prime }},\left( \rho ^{\prime },\eta ^{\prime }\right)
\right)
\end{equation*}%
target as being the set 
\begin{equation*}
\begin{array}{c}
\left\{ \left( \varphi ,\varphi _{0}\right) \in \mathbf{B}^{\mathbf{v}%
}\left( \left( F,\nu ,N\right) ,\left( F^{\prime },\nu ^{\prime },N^{\prime
}\right) \right) \right\}%
\end{array}%
\end{equation*}%
such that $\varphi _{0}\in Iso_{\mathbf{Man}}\left( N,N^{\prime }\right) $
and the $\mathbf{Mod}$-morphism $\Gamma \left( \varphi ,\varphi _{0}\right) $
is a $\mathbf{LieAlg}$-morphism of 
\begin{equation*}
\left( \Gamma \left( F,\nu ,N\right) ,+,\cdot ,\left[ ,\right] _{F,h}\right)
\end{equation*}%
source and 
\begin{equation*}
\left( \Gamma \left( F^{\prime },\nu ^{\prime },N^{\prime }\right) ,+,\cdot
, \left[ ,\right] _{F^{\prime },h^{\prime }}\right)
\end{equation*}%
target.

\textbf{Definition 2.8} For any $\mathbf{GLA}$-morphism $\left( \varphi
,\varphi _{0}\right) $ of 
\begin{equation*}
\left( \left( F,\nu ,N\right) ,\left[ ,\right] _{F,h},\left( \rho ,\eta
\right) \right)
\end{equation*}%
source and 
\begin{equation*}
\left( \left( F^{\prime },\nu ^{\prime },N^{\prime }\right) ,\left[ ,\right]
_{F^{\prime },h^{\prime }},\left( \rho ^{\prime },\eta ^{\prime }\right)
\right)
\end{equation*}%
target we define the application 
\begin{equation*}
\begin{array}{ccc}
\Lambda ^{q}\left( F^{\prime },\nu ^{\prime },N^{\prime }\right) & ^{%
\underrightarrow{\ \left( \varphi ,\varphi _{0}\right) ^{\ast }\ }} & 
\Lambda ^{q}\left( F,\nu ,N\right) \\ 
\omega ^{\prime } & \longmapsto & \left( \varphi ,\varphi _{0}\right) ^{\ast
}\omega ^{\prime }%
\end{array}%
,
\end{equation*}%
where 
\begin{equation*}
\begin{array}{c}
\left( \left( \varphi ,\varphi _{0}\right) ^{\ast }\omega ^{\prime }\right)
\left( z_{1},...,z_{q}\right) =\omega ^{\prime }\left( \Gamma \left( \varphi
,\varphi _{0}\right) \left( z_{1}\right) ,...,\Gamma \left( \varphi ,\varphi
_{0}\right) \left( z_{q}\right) \right) ,%
\end{array}%
\end{equation*}%
for any $z_{1},...,z_{q}\in \Gamma \left( F,\nu ,N\right) .$

\textit{\noindent Remark 2.5}\textbf{\ }It is remarked that the $\mathbf{B}^{%
\mathbf{v}}$-morphism $\left( Th\circ \rho ,h\circ \eta \right) $ is a $%
\mathbf{GLA}$-morphism~of 
\begin{equation*}
\left( \left( F,\nu ,N\right) ,\left[ ,\right] _{F,h},\left( \rho ,\eta
\right) \right)
\end{equation*}%
\emph{\ }source and 
\begin{equation*}
\left( \left( TN,\tau _{N},N\right) ,\left[ ,\right] _{TN,Id_{N}},\left(
Id_{TN},Id_{N}\right) \right)
\end{equation*}%
target.

Moreover, for any $\tilde{\imath}\in \overline{1,n}$, we obtain%
\begin{equation*}
\begin{array}{c}
\left( Th\circ \rho ,h\circ \eta \right) ^{\ast }\left( d\varkappa ^{\tilde{%
\imath}}\right) =d^{F}\varkappa ^{\tilde{\imath}},%
\end{array}%
\end{equation*}%
where $d$ is the exterior differentiation operator associated to the
exterior differential Lie $\mathcal{F}\left( N\right) $-algebra\emph{\ }%
\begin{equation*}
\begin{array}{c}
\left( \Lambda \left( TN,\tau _{N},N\right) ,+,\cdot ,\wedge \right) .%
\end{array}%
\end{equation*}

\textbf{Theorem 2.11 }\emph{If }$\left( \varphi ,\varphi _{0}\right) $\emph{%
\ is a morphism of }%
\begin{equation*}
\left( \left( F,\nu ,N\right) ,\left[ ,\right] _{F,h},\left( \rho ,\eta
\right) \right)
\end{equation*}%
\emph{\ source and }%
\begin{equation*}
\left( \left( F^{\prime },\nu ^{\prime },N^{\prime }\right) ,\left[ ,\right]
_{F^{\prime },h^{\prime }},\left( \rho ^{\prime },\eta ^{\prime }\right)
\right)
\end{equation*}%
\emph{target, then the following affirmations are satisfied:}\medskip

\noindent 1. \emph{For any }$\omega ^{\prime }\in \Lambda ^{q}\left(
F^{\prime },\nu ^{\prime },N^{\prime }\right) $\emph{\ and }$\theta ^{\prime
}\in \Lambda ^{r}\left( F^{\prime },\nu ^{\prime },N^{\prime }\right) $\emph{%
\ we obtain}%
\begin{equation*}
\begin{array}{c}
\left( \varphi ,\varphi _{0}\right) ^{\ast }\left( \omega ^{\prime }\wedge
\theta ^{\prime }\right) =\left( \varphi ,\varphi _{0}\right) ^{\ast }\omega
^{\prime }\wedge \left( \varphi ,\varphi _{0}\right) ^{\ast }\theta ^{\prime
}.%
\end{array}%
\leqno(2.18)
\end{equation*}

\noindent 2.\emph{\ For any }$z\in \Gamma \left( F,\nu ,N\right) $\emph{\
and }$\omega ^{\prime }\in \Lambda ^{q}\left( F^{\prime },\nu ^{\prime
},N^{\prime }\right) $\emph{\ we obtain}%
\begin{equation*}
\begin{array}{c}
i_{z}\left( \left( \varphi ,\varphi _{0}\right) ^{\ast }\omega ^{\prime
}\right) =\left( \varphi ,\varphi _{0}\right) ^{\ast }\left( i_{\Gamma
\left( \varphi ,\varphi _{0}\right) z}\omega ^{\prime }\right) .%
\end{array}%
\leqno(2.19)
\end{equation*}

\noindent 3. \emph{If }$N=N^{\prime }$ \emph{and }%
\begin{equation*}
\left( Th\circ \rho ,h\circ \eta \right) =\left( Th^{\prime }\circ \rho
^{\prime },h^{\prime }\circ \eta ^{\prime }\right) \circ \left( \varphi
,\varphi _{0}\right) ,
\end{equation*}%
\emph{then we obtain} 
\begin{equation*}
\begin{array}{c}
\left( \varphi ,\varphi _{0}\right) ^{\ast }\circ d^{F^{\prime }}=d^{F}\circ
\left( \varphi ,\varphi _{0}\right) ^{\ast }.%
\end{array}%
\leqno(2.20)
\end{equation*}

\emph{Proof.} 1. Let $\omega ^{\prime }\in \Lambda ^{q}\left( F^{\prime
},\nu ^{\prime },N^{\prime }\right) $ and $\theta ^{\prime }\in \Lambda
^{r}\left( F^{\prime },\nu ^{\prime },N^{\prime }\right) $ be arbitrary.
Since 
\begin{equation*}
\begin{array}{l}
\displaystyle\left( \varphi ,\varphi _{0}\right) ^{\ast }\left( \omega
^{\prime }\wedge \theta ^{\prime }\right) \left( z_{1},...,z_{q+r}\right)
=\left( \omega ^{\prime }\wedge \theta ^{\prime }\right) \left( \Gamma
\left( \varphi ,\varphi _{0}\right) z_{1},...,\Gamma \left( \varphi ,\varphi
_{0}\right) z_{q+r}\right) \vspace*{1mm} \\ 
\qquad\displaystyle=\frac{1}{\left( q+r\right) !}\underset{\sigma \in \Sigma
_{q+r}}{\tsum }sgn\left( \sigma \right) \cdot \omega ^{\prime }\left( \Gamma
\left( \varphi ,\varphi _{0}\right) z_{1},...,\Gamma \left( \varphi ,\varphi
_{0}\right) z_{q}\right) \vspace*{1mm} \\ 
\hfill\cdot \theta ^{\prime }\left( \Gamma \left( \varphi ,\varphi
_{0}\right) z_{q+1},...,\Gamma \left( \varphi ,\varphi _{0}\right)
z_{q+r}\right) \vspace*{1mm} \\ 
\qquad\displaystyle=\frac{1}{\left( q+r\right) !}\underset{\sigma \in \Sigma
_{q+r}}{\tsum }sgn\left( \sigma \right) \cdot \left( \varphi ,\varphi
_{0}\right) ^{\ast }\omega ^{\prime }\left( z_{1},...,z_{q}\right) \left(
\varphi ,\varphi _{0}\right) ^{\ast }\theta ^{\prime }\left(
z_{q+1},...,z_{q+r}\right) \vspace*{1mm} \\ 
\qquad\displaystyle=\left( \left( \varphi ,\varphi _{0}\right) ^{\ast
}\omega ^{\prime }\wedge \left( \varphi ,\varphi _{0}\right) ^{\ast }\theta
^{\prime }\right) \left( z_{1},...,z_{q+r}\right) ,%
\end{array}%
\end{equation*}%
for any $z_{1},...,z_{q+r}\in \Gamma \left( F,\nu ,N\right) $, it results
the conclusion of affirmation 1.\medskip

2. Let $z\in \Gamma \left( F,\nu ,N\right) $\emph{\ }and\emph{\ }$\omega
^{\prime }\in \Lambda ^{q}\left( F^{\prime },\nu ^{\prime },N^{\prime
}\right) $ be arbitrary. Since 
\begin{equation*}
\begin{array}{ll}
i_{z}\left( \left( \varphi ,\varphi _{0}\right) ^{\ast }\omega ^{\prime
}\right) \left( z_{2},...,z_{q}\right) & \displaystyle =\omega ^{\prime
}\left( \Gamma \left( \varphi ,\varphi _{0}\right) z,\Gamma \left( \varphi
,\varphi _{0}\right) z_{2},...,\Gamma \left( \varphi ,\varphi _{0}\right)
z_{q}\right) \vspace*{1mm} \\ 
& \displaystyle=i_{\Gamma \left( \varphi ,\varphi _{0}\right) z}\omega
^{\prime }\left( \Gamma \left( \varphi ,\varphi _{0}\right) z_{2},...,\Gamma
\left( \varphi ,\varphi _{0}\right) z_{q}\right) \vspace*{1mm} \\ 
& \displaystyle=\left( \varphi ,\varphi _{0}\right) ^{\ast }\left( i_{\Gamma
\left( \varphi ,\varphi _{0}\right) z}\omega ^{\prime }\right) \left(
z_{2},...,z_{q}\right) ,%
\end{array}%
\end{equation*}%
for any $z_{2},...,z_{q}\in \Gamma \left( F,\nu ,N\right) $, it results the
conclusion of affirmation 2$.$\smallskip

3. Let $\omega ^{\prime }\in \Lambda ^{q}\left( F^{\prime },\nu ^{\prime
},N^{\prime }\right) $ and $z_{0},...,z_{q}\in \Gamma \left( F,\nu ,N\right) 
$ be arbitrary. Since 
\begin{equation*}
\begin{array}{l}
\left( \left( \varphi ,\varphi _{0}\right) ^{\ast }d^{F^{\prime }}\omega
^{\prime }\right) \left( z_{0},...,z_{q}\right) =\left( d^{F^{\prime
}}\omega ^{\prime }\right) \left( \Gamma \left( \varphi ,\varphi _{0}\right)
z_{0},...,\Gamma \left( \varphi ,\varphi _{0}\right) z_{q}\right) \vspace*{%
1mm} \\ 
=\overset{q}{\underset{i=0}{\tsum }}\left( -1\right) ^{i}\Gamma \left(
Th^{\prime }\circ \rho ^{\prime },h^{\prime }\circ \eta ^{\prime }\right)
\left( \Gamma \left( \varphi ,\varphi _{0}\right) z_{i}\right) \vspace*{1mm}
\\ 
\hfill \cdot \omega ^{\prime }\left( \left( \Gamma \left( \varphi ,\varphi
_{0}\right) z_{0},\Gamma \left( \varphi ,\varphi _{0}\right) z_{1},...,%
\widehat{\Gamma \left( \varphi ,\varphi _{0}\right) z_{i}},...,\Gamma \left(
\varphi ,\varphi _{0}\right) z_{q}\right) \right) \vspace*{1mm} \\ 
+\underset{0\leq i<j}{\tsum }\left( -1\right) ^{i+j}\cdot \omega ^{\prime
}\left( \Gamma \left( \varphi ,\varphi _{0}\right) \left[ z_{i},z_{j}\right]
_{F,h},\Gamma \left( \varphi ,\varphi _{0}\right) z_{0},\Gamma \left(
\varphi ,\varphi _{0}\right) z_{1},...,\right. \vspace*{1mm} \\ 
\hfill \left. \cdot \widehat{\Gamma \left( \varphi ,\varphi _{0}\right) z_{i}%
},...,\widehat{\Gamma \left( \varphi ,\varphi _{0}\right) z_{j}},...,\Gamma
\left( \varphi ,\varphi _{0}\right) z_{q}\right)%
\end{array}%
\end{equation*}%
and 
\begin{equation*}
\begin{array}{l}
d^{F}\left( \left( \varphi ,\varphi _{0}\right) ^{\ast }\omega ^{\prime
}\right) \left( z_{0},...,z_{q}\right) \vspace*{1mm} \\ 
=\overset{q}{\underset{i=0}{\tsum }}\left( -1\right) ^{i}\Gamma \left(
Th\circ \rho ,h\circ \eta \right) \left( z_{i}\right) \cdot \left( \left(
\varphi ,\varphi _{0}\right) ^{\ast }\omega ^{\prime }\right) \left(
z_{0},...,\widehat{z_{i}},...,z_{q}\right) \vspace*{1mm} \\ 
+\underset{0\leq i<j}{\tsum }\left( -1\right) ^{i+j}\cdot \left( \left(
\varphi ,\varphi _{0}\right) ^{\ast }\omega ^{\prime }\right) \left( \left[
z_{i},z_{j}\right] _{F,h},z_{0},...,\widehat{z_{i}},...,\widehat{z_{j}}%
,...,z_{q}\right) \vspace*{1mm} \\ 
=\overset{q}{\underset{i=0}{\tsum }}\left( -1\right) ^{i}\Gamma \left(
Th\circ \rho ,h\circ \eta \right) \left( z_{i}\right) \cdot \omega ^{\prime
}\left( \Gamma \left( \varphi ,\varphi _{0}\right) z_{0},...,\widehat{\Gamma
\left( \varphi ,\varphi _{0}\right) z_{i}},...,\Gamma \left( \varphi
,\varphi _{0}\right) z_{q}\right) \vspace*{1mm} \\ 
+\underset{0\leq i<j}{\tsum }\left( -1\right) ^{i+j}\cdot \omega ^{\prime
}\left( \Gamma \left( \varphi ,\varphi _{0}\right) \left[ z_{i},z_{j}\right]
_{F,h},\Gamma \left( \varphi ,\varphi _{0}\right) z_{0},\Gamma \left(
\varphi ,\varphi _{0}\right) z_{1},...,\right. \vspace*{1mm} \\ 
\left. \widehat{\Gamma \left( \varphi ,\varphi _{0}\right) z_{i}},...,%
\widehat{\Gamma \left( \varphi ,\varphi _{0}\right) z_{j}},...,\Gamma \left(
\varphi ,\varphi _{0}\right) z_{q}\right)%
\end{array}%
\end{equation*}%
it results the conclusion of affirmation 3. \hfill \emph{q.e.d.}

\textbf{Definition 2.9} For any $q\in \overline{1,n}$ we define%
\begin{equation*}
\begin{array}{c}
\mathcal{Z}^{q}\left( F,\nu ,N\right) =\left\{ \omega \in \Lambda ^{q}\left(
F,\nu ,N\right) :d^{F}\omega =0\right\} ,%
\end{array}%
\end{equation*}%
the set of \emph{closed differential exterior }$q$\emph{-forms} and 
\begin{equation*}
\begin{array}{c}
\mathcal{B}^{q}\left( F,\nu ,N\right) =\left\{ \omega \in \Lambda ^{q}\left(
F,\nu ,N\right) :\exists \eta \in \Lambda ^{q-1}\left( F,\nu ,N\right)
~|~d^{F}\eta =\omega \right\} ,%
\end{array}%
\end{equation*}%
the set of \emph{exact differential exterior }$q$\emph{-forms}.

\bigskip

\section{Torsion and curvature forms. Identities of Cartan and Bianchi type}

\ \ \ 

Using the theory of linear connections of Eresmann type presented in $\left[
1\right] $ for the diagram: 
\begin{equation*}
\begin{array}{rcl}
E &  & \left( F,\left[ ,\right] _{F,h},\left( \rho ,Id_{N}\right) \right) \\ 
\pi \downarrow &  & ~\downarrow \nu \\ 
M & ^{\underrightarrow{~\ \ \ \ h~\ \ \ \ }} & ~\ N%
\end{array}%
,\leqno(3.1)
\end{equation*}%
where $\left( E,\pi ,M\right) \in \left\vert \mathbf{B}^{\mathbf{v}%
}\right\vert $ and $\left( \left( F,\nu ,N\right) ,\left[ ,\right]
_{F,h},\left( \rho ,Id_{N}\right) \right) \in \left\vert \mathbf{GLA}%
\right\vert ,$ we obtain a linear $\rho $-connection $\rho \Gamma $ for the
vector bundle $\left( E,\pi ,M\right) $ by components $\rho \Gamma _{b\alpha
}^{a}.$

Using the components of this linear $\rho $-connection, we obtain a linear $%
\rho $-connection $\rho \dot{\Gamma}$ for the vector bundle $\left( E,\pi
,M\right) $ given by the diagram: 
\begin{equation*}
\begin{array}{ccl}
~\ \ \ E &  & \left( h^{\ast }F,\left[ ,\right] _{h^{\ast }F},\left( \overset%
{h^{\ast }F}{\rho },Id_{M}\right) \right) \\ 
\pi \downarrow &  & ~\ \ \downarrow h^{\ast }\nu \\ 
~\ \ \ M & ^{\underrightarrow{~\ \ \ \ Id_{M}~\ \ }} & ~\ \ \ M%
\end{array}%
\leqno(3.2)
\end{equation*}

If $\left( E,\pi ,M\right) =\left( F,\nu ,N\right) ,$ then, using the
components of the same linear $\rho $-connection $\rho \Gamma ,$ we can
consider a linear $\rho $-connection $\rho \ddot{\Gamma}$ for the vector
bundle $\left( h^{\ast }E,h^{\ast }\pi ,M\right) $ given by the diagram: 
\begin{equation*}
\begin{array}{ccl}
~\ \ ~\ \ \ \ h^{\ast }E &  & \left( h^{\ast }E,\left[ ,\right] _{h^{\ast
}E},\left( \overset{h^{\ast }E}{\rho },Id_{M}\right) \right) \\ 
h^{\ast }\pi \downarrow &  & ~\ \ \downarrow h^{\ast }\pi \\ 
~\ \ \ \ \ M & ^{\underrightarrow{~\ \ \ \ Id_{M}~\ \ }} & ~\ \ M%
\end{array}%
\leqno(3.3)
\end{equation*}

\textbf{Definitiona 3.1 }For any $a,b\in \overline{1,n}$ we define the
connection form $\Omega _{b}^{a}=\rho \Gamma _{bc}^{a}S^{c}.$\bigskip
\noindent\ In the particular case of Lie algebroid, $h=Id_{M},$ we obtain
the connection form $\omega _{b}^{a}=\rho \Gamma _{bc}^{a}s^{c}.$ In
addition, if $h=Id_{M},$ then we obtain the connection form $\omega
_{j}^{i}=\Gamma _{jk}^{i}dx^{k}.$

\textbf{Definition 3.2 }If $\left( E,\pi ,M\right) =\left( F,\nu ,N\right) $%
, then the application 
\begin{equation*}
\begin{array}{ccc}
\Gamma \left( h^{\ast }E,h^{\ast }\pi ,M\right) ^{2} & ^{\underrightarrow{\
\left( \rho ,h\right) \mathbb{T}\ }} & \Gamma \left( h^{\ast }E,h^{\ast }\pi
,M\right) \\ 
\left( U,V\right) & \longrightarrow & \rho \mathbb{T}\left( U,V\right)%
\end{array}%
\leqno(3.4)
\end{equation*}%
defined by: 
\begin{equation*}
\begin{array}{c}
\left( \rho ,h\right) \mathbb{T}\left( U,V\right) =\rho \ddot{D}_{U}V-\rho 
\ddot{D}_{V}U-\left[ U,V\right] _{h^{\ast }E},\,%
\end{array}%
\leqno(5.5)
\end{equation*}%
for any $U,V\in \Gamma \left( h^{\ast }E,h^{\ast }\pi ,M\right) ,$ is the $%
\left( \rho ,h\right) $\emph{-}torsion associated to linear $\rho $%
-connection $\rho \Gamma .$

If\emph{\ }$\left( \rho ,h\right) \mathbb{T}\left( S_{a},S_{b}\right) 
\overset{put}{=}\left( \rho ,h\right) \mathbb{T}_{~ab}^{c}S_{c},$\emph{\ }%
then the vector valued $2$-form%
\begin{equation*}
\begin{array}{c}
\left( \rho ,h\right) \mathbb{T}=\left( \left( \rho ,h\right) \mathbb{T}%
_{~ab}^{c}S_{c}\right) S^{a}\wedge S^{b}%
\end{array}%
\leqno(3.6)
\end{equation*}%
will be called the \emph{vector valued form of }$\left( \rho ,h\right) $%
\emph{-torsion }$\left( \rho ,h\right) \mathbb{T}$\emph{.}

In the particular case of Lie algebroids, $h=Id_{M}$, then the vector valued 
$2$-form becomes: 
\begin{equation*}
\begin{array}{c}
\rho \mathbb{T}=\left( \rho \mathbb{T}_{~ab}^{c}s_{c}\right) s^{a}\wedge
s^{b}.%
\end{array}%
\leqno(3.6)^{\prime }
\end{equation*}

In the classical case, $\rho =Id_{TM}$, then the vector valued $2$-form $%
\left( 3.6^{\prime }\right) $ becomes: 
\begin{equation*}
\begin{array}{c}
\mathbb{T}=\left( \mathbb{T}_{~jk}^{i}\frac{\partial }{\partial x^{i}}%
\right) dx^{j}\wedge dx^{k}.%
\end{array}%
\leqno(3.6)^{\prime \prime }
\end{equation*}

\textbf{Definition 3.3 }For each $c\in \overline{1,n}$ we obtain the \emph{%
scalar }$2$\emph{-form of }$\left( \rho ,h\right) $\emph{-torsion }$\left(
\rho ,h\right) \mathbb{T}$ 
\begin{equation*}
\begin{array}{c}
\left( \rho ,h\right) \mathbb{T}^{c}=\left( \rho ,h\right) \mathbb{T}%
_{~ab}^{c}S^{a}\wedge S^{b}.%
\end{array}%
\leqno(3.7)
\end{equation*}

In the particular case of Lie algebroids, $h=Id_{M}$, then the scalar $2$%
-form $(3.7)$ becomes: 
\begin{equation*}
\begin{array}{c}
\rho \mathbb{T}^{c}=\rho \mathbb{T}_{~ab}^{c}s^{a}\wedge s^{b}.%
\end{array}%
\leqno(3.7)^{\prime \prime }
\end{equation*}

In the classical case, $\rho =Id_{TM}$, then the scalar $2$-form $\left(
3.7^{\prime }\right) $ becomes: 
\begin{equation*}
\begin{array}{c}
\mathbb{T}^{i}=\mathbb{T}_{~jk}^{i}dx^{j}\wedge dx^{k}.%
\end{array}%
\leqno(3.7)^{\prime \prime }
\end{equation*}

\textbf{Definition 3.4 }The application 
\begin{equation*}
\begin{array}{ccl}
(\Gamma \left( h^{\ast }F,h^{\ast }\nu ,M\right) ^{2}{\times }\Gamma (E,\pi
,M) & ^{\underrightarrow{\ \left( \rho ,h\right) \mathbb{R}\ }} & \Gamma
(E,\pi ,M) \\ 
((Z,V),u) & \longrightarrow & \rho \mathbb{R}(Z,V)u%
\end{array}%
\leqno(3.8)
\end{equation*}%
defined by 
\begin{equation*}
\left( \rho ,h\right) \mathbb{R}\left( Z,V\right) u=\rho \dot{D}_{Z}\left(
\rho \dot{D}_{V}u\right) -\rho \dot{D}_{V}\left( \rho \dot{D}_{Z}u\right)
-\rho \dot{D}_{\left[ Z,V\right] _{h^{\ast }F}}u,\,\leqno(3.9)
\end{equation*}%
for any $Z,V\in \Gamma \left( h^{\ast }F,h^{\ast }\nu ,M\right) ,~u\in
\Gamma \left( E,\pi ,M\right) ,$ is called $\left( \rho ,h\right) $%
-curvature associated to linear $\rho $-connection $\rho \Gamma .$

If 
\begin{equation*}
\left( \rho ,h\right) \mathbb{R}\left( T_{\beta },T_{\alpha }\right) s_{b}%
\overset{put}{=}\left( \rho ,h\right) \mathbb{R}_{b~\alpha \beta }^{a}s_{a},
\end{equation*}%
then the vector mixed form 
\begin{equation*}
\begin{array}{c}
\left( \rho ,h\right) \mathbb{R=}\left( \left( \left( \rho ,h\right) \mathbb{%
R}_{b~\alpha \beta }^{a}s_{a}\right) T^{\alpha }\wedge T^{\beta }\right)
s^{b}%
\end{array}%
\leqno(3.10)
\end{equation*}%
will be called the \emph{vector valued form of }$\left( \rho ,h\right) $%
\emph{-curvature }$\left( \rho ,h\right) \mathbb{R}$\emph{.}

In the particular case of Lie algebroids, $h=Id_{M}$, then the vector mixed
form $\left( 3.10\right) $ becomes: 
\begin{equation*}
\begin{array}{c}
\rho \mathbb{R=}\left( \left( \rho \mathbb{R}_{b~\alpha \beta
}^{a}s_{a}\right) t^{\alpha }\wedge t^{\beta }\right) s^{b}%
\end{array}%
\leqno(3.10)^{\prime }
\end{equation*}

In the classical case, $h=Id_{M}$, then the vector mixed form $\left(
3.10\right) ^{\prime }$ becomes:%
\begin{equation*}
\begin{array}{c}
\mathbb{R=}\left( \left( \mathbb{R}_{b~hk}^{a}s_{a}\right) dx^{h}\wedge
dx^{k}\right) s^{b}.%
\end{array}%
\leqno(3.10)^{\prime \prime }
\end{equation*}

\textbf{Definition 3.5 }For each $a,b\in \overline{1,n}$ we obtain the \emph{%
scalar }$2$\emph{-form of }$\left( \rho ,h\right) $\emph{-curvature }$\left(
\rho ,h\right) \mathbb{R}$ 
\begin{equation*}
\begin{array}{c}
\left( \rho ,h\right) \mathbb{R}_{b}^{a}=\left( \rho ,h\right) \mathbb{R}%
_{b~\alpha \beta }^{a}T^{\alpha }\wedge T^{\beta }.%
\end{array}%
\leqno(3.11)
\end{equation*}

In the particular case of Lie algebroids, $h=Id_{M}$, the \emph{scalar }$2$%
\emph{-form }$(3.11)\mathbb{\ }$becomes 
\begin{equation*}
\begin{array}{c}
\rho \mathbb{R}_{b}^{a}=\rho \mathbb{R}_{b~\alpha \beta }^{a}t^{\alpha
}\wedge t^{\beta }.%
\end{array}%
\leqno(3.11)^{\prime }
\end{equation*}

In the classical case, $h=Id_{M}$, the \emph{scalar }$2$\emph{-form }$%
(3.11)^{\prime }\mathbb{\ }$becomes: 
\begin{equation*}
\begin{array}{c}
\mathbb{R}_{b}^{a}=\mathbb{R}_{b~hk}^{a}dx^{h}\wedge dx^{k}.%
\end{array}%
\leqno(3.11)^{\prime \prime }
\end{equation*}

\textbf{Theorem 3.1 }\emph{The identities }%
\begin{equation*}
\begin{array}{c}
\left( \rho ,h\right) \mathbb{T}^{a}=d^{h^{\ast }F}S^{a}+\Omega
_{b}^{a}\wedge S^{b},%
\end{array}%
\leqno(C_{1})
\end{equation*}%
\emph{and }%
\begin{equation*}
\begin{array}{c}
\left( \rho ,h\right) \mathbb{R}_{b}^{a}=d^{h^{\ast }F}\Omega
_{b}^{a}+\Omega _{c}^{a}\wedge \Omega _{b}^{c}%
\end{array}%
\leqno(C_{2})
\end{equation*}%
\emph{hold good. These will be called the first respectively the second
identity of Cartan type.}

\emph{In the particular case of Lie algebroids, }$h=Id_{M}$, \emph{then the
identities }$(C_{1})$\emph{\ and }$(C_{2})$ \emph{become}%
\begin{equation*}
\begin{array}{c}
\rho \mathbb{T}^{a}=d^{F}s^{a}+\omega _{b}^{a}\wedge s^{b},%
\end{array}%
\leqno(C_{1}^{\prime })
\end{equation*}%
\emph{and }%
\begin{equation*}
\begin{array}{c}
\rho \mathbb{R}_{b}^{a}=d^{F}\omega _{b}^{a}+\omega _{c}^{a}\wedge \omega
_{b}^{c}%
\end{array}%
\leqno(C_{2}^{\prime })
\end{equation*}%
\emph{respectively.}

\emph{In the classical case,}$\rho =Id_{TM}$\emph{, then the identities }$%
\left( C_{1}^{\prime }\right) $\emph{\ and }$\left( C_{2}^{\prime }\right)
\, $\emph{become:}%
\begin{equation*}
\begin{array}{c}
\mathbb{T}^{i}=ddx^{i}+\omega _{j}^{i}\wedge dx^{j}=\omega _{j}^{i}\wedge
dx^{j}%
\end{array}%
\leqno(C_{1}^{\prime \prime })
\end{equation*}%
\emph{and }%
\begin{equation*}
\begin{array}{c}
\mathbb{R}_{j}^{i}=d\omega _{j}^{i}+\omega _{h}^{i}\wedge \omega _{j}^{h},%
\end{array}%
\leqno(C_{2}^{\prime \prime })
\end{equation*}%
\emph{\ respectively. }\hfill

\emph{Proof.} To prove the first identity we consider that $\left( E,\pi
,M\right) =\left( F,\nu ,M\right) .$ Since 
\begin{equation*}
\begin{array}{l}
d^{h^{\ast }F}S^{a}(U,V)S_{a}=((\Gamma (\overset{h^{\ast }F}{\rho }%
,Id_{M})U)S^{a}(V)\vspace*{1mm} \\ 
\qquad -(\Gamma (\overset{h^{\ast }F}{\rho }%
,Id_{M})V)(S^{a}(U))-S^{a}([U,V]_{h^{\ast }F}))S_{a}\vspace*{1mm} \\ 
\qquad =(\Gamma (\overset{h^{\ast }F}{\rho },Id_{M})U)(V^{a})-(\Gamma (%
\overset{h^{\ast }F}{\rho },Id_{M})V)(U^{a})-S^{a}([U,V]_{h^{\ast }F})S_{a}%
\vspace*{1mm} \\ 
\qquad =\rho \ddot{D}_{U}V-V^{b}\rho \ddot{D}_{U}S_{b}-\rho \ddot{D}%
_{V}U-U^{b}\rho \ddot{D}_{V}S_{b}-[U,V]_{h^{\ast }F}\vspace*{1mm} \\ 
\qquad =\left( \rho ,h\right) \mathbb{T}(U,V)-(\rho \Gamma
_{bc}^{a}V^{b}U^{c}-\rho \Gamma _{bc}^{a}U^{b}V^{c})S_{a}\vspace*{1mm} \\ 
\qquad =(\left( \rho ,h\right) \mathbb{T}^{a}(U,V)-\Omega _{b}^{a}\wedge
S^{b}(U,V))S_{a},%
\end{array}%
\end{equation*}%
it results the first identity.

To prove the second identity, we consider that $\left( E,\pi ,M\right) \neq
\left( F,\nu ,M\right) .$ Since 
\begin{equation*}
\begin{array}{l}
\left( \rho ,h\right) \mathbb{R}_{b}^{a}\left( Z,W\right) s_{a}=\left( \rho
,h\right) \mathbb{R}\left( \left( W,Z\right) ,s_{b}\right) \vspace*{1mm} \\ 
\qquad =\rho \dot{D}_{Z}\left( \rho \dot{D}_{W}s_{b}\right) -\rho \dot{D}%
_{W}\left( \rho \dot{D}_{Z}s_{b}\right) -\rho \dot{D}_{\left[ Z,W\right]
_{h^{\ast }F}}s_{b}\vspace*{1mm} \\ 
\qquad =\rho \dot{D}_{Z}\left( \Omega _{b}^{a}\left( W\right) s_{a}\right)
-\rho \dot{D}_{W}\left( \Omega _{b}^{a}\left( Z\right) s_{a}\right) -\Omega
_{b}^{a}\left( \left[ Z,W\right] _{h^{\ast }F}\right) s_{a}\vspace*{1mm} \\ 
\qquad +\left( \Omega _{c}^{a}\left( Z\right) \Omega _{b}^{c}\left( W\right)
-\Omega _{c}^{a}\left( W\right) \Omega _{b}^{c}\left( Z\right) \right) s_{a}%
\vspace*{1mm} \\ 
\qquad =\left( d^{h^{\ast }F}\Omega _{b}^{a}\left( Z,W\right) +\Omega
_{c}^{a}\wedge \Omega _{b}^{c}\left( Z,W\right) \right) s_{a}%
\end{array}%
\end{equation*}%
it results the second identity.\hfill \hfill \emph{q.e.d.}

\textbf{Theorem 3.2 }\emph{The identities }%
\begin{equation*}
\begin{array}{c}
d^{h^{\ast }F}\left( \rho ,h\right) \mathbb{T}^{a}=\left( \rho ,h\right) 
\mathbb{R}_{b}^{a}\wedge S^{b}-\Omega _{c}^{a}\wedge \left( \rho ,h\right) 
\mathbb{T}^{c}%
\end{array}%
\leqno(B_{1})
\end{equation*}%
\emph{and }%
\begin{equation*}
\begin{array}{c}
d^{h^{\ast }F}\left( \rho ,h\right) \mathbb{R}_{b}^{a}=\left( \rho ,h\right) 
\mathbb{R}_{c}^{a}\wedge \Omega _{b}^{c}-\Omega _{c}^{a}\wedge \left( \rho
,h\right) \mathbb{R}_{b}^{c},%
\end{array}%
\leqno(B_{2})
\end{equation*}%
\emph{hold good.} \emph{We will called these the first respectively the
second identity of Bianchi type.}

\emph{If the }$\left( \rho ,h\right) $\emph{-torsion is null, then the first
identity of Bianchi type becomes: }%
\begin{equation*}
\begin{array}{c}
\left( \rho ,h\right) \mathbb{R}_{b}^{a}\wedge s^{b}=0.%
\end{array}%
\leqno(\tilde{B}_{1})
\end{equation*}

\emph{In the particularcase of Lie algebroids, }$h=Id_{M},$\emph{\ then the
identities }$(B_{1})$\emph{\ and }$(B_{2})$ \emph{become }%
\begin{equation*}
\begin{array}{c}
d^{F}\rho \mathbb{T}^{a}=\rho \mathbb{R}_{b}^{a}\wedge s^{b}-\omega
_{c}^{a}\wedge \rho \mathbb{T}^{c}%
\end{array}%
\leqno(B_{1}^{\prime })
\end{equation*}%
\emph{and }%
\begin{equation*}
\begin{array}{c}
d^{F}\rho \mathbb{R}_{b}^{a}=\rho \mathbb{R}_{c}^{a}\wedge \omega
_{b}^{c}-\omega _{c}^{a}\wedge \rho \mathbb{R}_{b}^{c},%
\end{array}%
\leqno(B_{2}^{\prime })
\end{equation*}%
\emph{\ respectively.}

\emph{In the classical case, }$\rho =Id_{TM}$\emph{, then the identities }$%
\left( B_{1}^{\prime }\right) $\emph{\ and }$\left( B_{2}^{\prime }\right)
\, $\emph{become:} 
\begin{equation*}
\begin{array}{c}
d\mathbb{T}^{i}=\mathbb{R}_{j}^{i}\wedge dx^{j}-\omega _{k}^{i}\wedge 
\mathbb{T}^{k}%
\end{array}%
\leqno(B_{1}^{\prime \prime })
\end{equation*}%
\emph{and }%
\begin{equation*}
\begin{array}{c}
d\mathbb{R}_{j}^{i}=\mathbb{R}_{h}^{i}\wedge \omega _{j}^{h}-\omega
_{h}^{i}\wedge \mathbb{R}_{j}^{h},%
\end{array}%
\leqno(B_{2}^{\prime \prime })
\end{equation*}%
\emph{\ respectively. }

\textit{Proof.} We consider $\left( E,\pi ,M\right) =\left( F,\nu ,M\right)
. $ Using the first identity of Cartan type and the equality $d^{h^{\ast
}F}\circ d^{h^{\ast }F}=0,$ we obtain: 
\begin{equation*}
\begin{array}{c}
d^{h^{\ast }F}\left( \rho ,h\right) \mathbb{T}^{a}=d^{h^{\ast }F}\Omega
_{b}^{a}\wedge S^{b}-\Omega _{c}^{a}\wedge d^{h^{\ast }F}S^{c}.%
\end{array}%
\end{equation*}

Using the second identity of Cartan type and the previous identity, we
obtain: 
\begin{equation*}
d^{h^{\ast }F}\left( \rho ,h\right) \mathbb{T}^{a}=\left( \left( \rho
,h\right) \mathbb{R}_{b}^{a}-\Omega _{c}^{a}\wedge \Omega _{b}^{c}\right)
\wedge S^{b}-\Omega _{c}^{a}\wedge \left( \left( \rho ,h\right) \mathbb{T}%
^{c}-\Omega _{b}^{c}\wedge S^{b}\right) .
\end{equation*}

After some calculations, we obtain the first identity of Bianchi type.

Using the second identity of Cartan type and the equality $d^{h^{\ast
}F}\circ d^{h^{\ast }F}=0,$ we obtain: 
\begin{equation*}
\begin{array}{c}
d^{h^{\ast }F}\Omega _{c}^{a}\wedge \Omega _{b}^{c}-\Omega _{c}^{a}\wedge
d^{h^{\ast }F}\Omega _{b}^{c}=d^{h^{\ast }F}\left( \rho ,h\right) \mathbb{R}%
_{b}^{a}.%
\end{array}%
\end{equation*}

Using the second of Cartan type and the previous identity, we obtain:%
\begin{equation*}
\begin{array}{cc}
d^{h^{\ast }F}\left( \rho ,h\right) \mathbb{R}_{b}^{a}=\left( \left( \rho
,h\right) \mathbb{R}_{c}^{a}-\Omega _{e}^{a}\wedge \Omega _{c}^{e}\right)
\wedge \Omega _{b}^{c}-\Omega _{c}^{a}\wedge \left( \left( \rho ,h\right) 
\mathbb{R}_{b}^{c}-\Omega _{e}^{c}\wedge \Omega _{b}^{e}\right) . & 
\end{array}%
\end{equation*}

After some calculations, we obtain the second identity of Bianchi
type.\hfill \hfill \emph{q.e.d.}

\bigskip

\section{Interior and exterior differential systems}

\ \ \ \ 

Let $\left( \left( F,\nu ,N\right) ,\left[ ,\right] _{F,h},\left( \rho ,\eta
\right) \right) $ be an object of the category $\mathbf{GLA}$.

Let $\mathcal{AF}_{F}$ be a vector fibred $\left( n+p\right) $-atlas for the
vector bundle $\left( F,\nu ,N\right) $ and let $\mathcal{AF}_{TM}$ be a
vector fibred $\left( m+m\right) $-atlas for the vector bundle $\left(
TM,\tau _{M},M\right) $.

Let $\left( h^{\ast }F,h^{\ast }\nu ,M\right) $ be the pull-back vector
bundle through $h.$

If $\left( U,\xi _{U}\right) \in \mathcal{AF}_{TM}$ and $\left(
V,s_{V}\right) \in \mathcal{AF}_{F}$ such that $U\cap h^{-1}\left( V\right)
\neq \phi $, then we define the application%
\begin{equation*}
\begin{array}{ccc}
h^{\ast }\nu ^{-1}(U{\cap }h^{-1}(V))) & {}^{\underrightarrow{\bar{s}_{U{%
\cap }h^{-1}(V)}}} & \left( U{\cap }h^{-1}(V)\right) {\times }\mathbb{R}^{p}
\\ 
\left( \varkappa ,z\left( h\left( \varkappa \right) \right) \right) & 
\longmapsto & \left( \varkappa ,t_{V,h\left( \varkappa \right) }^{-1}z\left(
h\left( \varkappa \right) \right) \right) .%
\end{array}%
\end{equation*}

\textbf{Proposition 4.1 }\emph{The set}%
\begin{equation*}
\begin{array}{c}
\overline{\mathcal{AF}}_{F}\overset{put}{=}\underset{U\cap h^{-1}\left(
V\right) \neq \phi }{\underset{\left( U,\xi _{U}\right) \in \mathcal{AF}%
_{TM},~\left( V,s_{V}\right) \in \mathcal{AF}_{F}}{\tbigcup }}\left\{ \left(
U\cap h^{-1}\left( V\right) ,\bar{s}_{U{\cap }h^{-1}(V)}\right) \right\}%
\end{array}%
\end{equation*}%
\emph{is a vector fibred }$m+p$\emph{-atlas for the vector bundle }$\left(
h^{\ast }F,h^{\ast }\nu ,M\right) .$

\emph{If }$z=z^{\alpha }t_{\alpha }\in \Gamma \left( F,\nu ,N\right) ,$ 
\emph{then\ we obtain the section }%
\begin{equation*}
\begin{array}{c}
Z=\left( z^{\alpha }\circ h\right) T_{\alpha }\in \Gamma \left( h^{\ast
}F,h^{\ast }\nu ,M\right)%
\end{array}%
\end{equation*}%
\emph{such that }$Z\left( x\right) =z\left( h\left( x\right) \right) ,$ 
\emph{for any }$x\in U\cap h^{-1}\left( V\right) .$\bigskip

\textbf{Theorem 4.1 }\emph{Let} $\Big({\overset{h^{\ast }F}{\rho }},Id_{M}%
\Big)$ \emph{be the }$\mathbf{B}^{\mathbf{v}}$\emph{-morphism of }$\left(
h^{\ast }F,h^{\ast }\nu ,M\right) $\ \emph{source and} $\left( TM,\tau
_{M},M\right) $\ \emph{target, where}%
\begin{equation*}
\begin{array}{rcl}
h^{\ast }F & ^{\underrightarrow{\overset{h^{\ast }F}{\rho }}} & TM \\ 
\displaystyle Z^{\alpha }T_{\alpha }\left( x\right) & \longmapsto & %
\displaystyle\left( Z^{\alpha }\cdot \rho _{\alpha }^{i}\circ h\right) \frac{%
\partial }{\partial x^{i}}\left( x\right)%
\end{array}%
\leqno(4.1)
\end{equation*}

\emph{Using the operation} 
\begin{equation*}
\begin{array}{ccc}
\Gamma \left( h^{\ast }F,h^{\ast }\nu ,M\right) \times \Gamma \left( h^{\ast
}F,h^{\ast }\nu ,M\right) & ^{\underrightarrow{~\ \ \left[ ,\right]
_{h^{\ast }F}~\ \ }} & \Gamma \left( h^{\ast }F,h^{\ast }\nu ,M\right)%
\end{array}%
\end{equation*}%
\emph{defined by}%
\begin{equation*}
\begin{array}{ll}
\left[ T_{\alpha },T_{\beta }\right] _{h^{\ast }F} & =\left( L_{\alpha \beta
}^{\gamma }\circ h\right) T_{\gamma },\vspace*{1mm} \\ 
\left[ T_{\alpha },fT_{\beta }\right] _{h^{\ast }F} & \displaystyle=f\left(
L_{\alpha \beta }^{\gamma }\circ h\right) T_{\gamma }+\left( \rho _{\alpha
}^{i}\circ h\right) \frac{\partial f}{\partial x^{i}}T_{\beta },\vspace*{1mm}
\\ 
\left[ fT_{\alpha },T_{\beta }\right] _{h^{\ast }F} & =-\left[ T_{\beta
},fT_{\alpha }\right] _{h^{\ast }F},%
\end{array}%
\leqno(4.2)
\end{equation*}%
\emph{for any} $f\in \mathcal{F}\left( M\right) ,$ \emph{it results that} 
\begin{equation*}
\begin{array}{c}
\left( \left( h^{\ast }F,h^{\ast }\nu ,M\right) ,\left[ ,\right] _{h^{\ast
}F},\left( \overset{h^{\ast }F}{\rho },Id_{M}\right) \right)%
\end{array}%
\end{equation*}%
is a Lie algebroid which is called \emph{the pull-back Lie algebroid of the
generalized Lie algebroid }$\left( \left( F,\nu ,N\right) ,\left[ ,\right]
_{F,h},\left( \rho ,\eta \right) \right) .$

\textbf{Definition 4.1 }Any vector subbundle $\left( E,\pi ,M\right) $ of
the pull-back vector bundle $\left( h^{\ast }F,h^{\ast }\nu ,M\right) $ will
be called \emph{interior differential system (IDS) of the generalized Lie
algebroid }%
\begin{equation*}
\left( \left( F,\nu ,N\right) ,\left[ ,\right] _{F,h},\left( \rho ,\eta
\right) \right) .
\end{equation*}

In particular, if $h=Id_{N}=\eta $, then we obtain the definition of \emph{%
IDS} of a Lie algebroid. (see $\left[ {2}\right] $)

\emph{Remark 4.1 }If $\left( E,\pi ,M\right) $ is an IDS of the generalized
Lie algebroid 
\begin{equation*}
\left( \left( F,\nu ,N\right) ,\left[ ,\right] _{F,h},\left( \rho ,\eta
\right) \right) ,
\end{equation*}%
then we obtain a vector subbundle $\left( E^{0},\pi ^{0},M\right) $ of the
vector bundle $\left( \overset{\ast }{h^{\ast }F},\overset{\ast }{h^{\ast
}\nu },M\right) $ such that 
\begin{equation*}
\Gamma \left( E^{0},\pi ^{0},M\right) \overset{put}{=}\left\{ \Omega \in
\Gamma \left( \overset{\ast }{h^{\ast }F},\overset{\ast }{h^{\ast }\nu }%
,M\right) :\Omega \left( S\right) =0,~\forall S\in \Gamma \left( E,\pi
,M\right) \right\} .
\end{equation*}

The vector subbundle $\left( E^{0},\pi ^{0},M\right) $ will be called \emph{%
the annihilator vector subbundle of the IDS }$\left( E,\pi ,M\right) .$

\textbf{Proposition 4.2 }\emph{If }$\left( E,\pi ,M\right) $\emph{\ is an
IDS of the generalized Lie algebroid }%
\begin{equation*}
\left( \left( F,\nu ,N\right) ,\left[ ,\right] _{F,h},\left( \rho ,\eta
\right) \right)
\end{equation*}%
\emph{such that }$\Gamma \left( E,\pi ,M\right) =\left\langle
S_{1},...,S_{r}\right\rangle $\emph{, then it exists }$\Theta
^{r+1},...,\Theta ^{p}\in \Gamma \left( \overset{\ast }{h^{\ast }F},\overset{%
\ast }{h^{\ast }\nu },M\right) $\emph{\ linearly independent such that }$%
\Gamma \left( E^{0},\pi ^{0},M\right) =\left\langle \Theta ^{r+1},...,\Theta
^{p}\right\rangle .$

\textbf{Definition 4.2 }The \emph{IDS} $\left( E,\pi ,M\right) $ of the
generalized Lie algebroid 
\begin{equation*}
\left( \left( F,\nu ,N\right) ,\left[ ,\right] _{F,h},\left( \rho ,\eta
\right) \right)
\end{equation*}%
will be called \emph{involutive} if $\left[ S,T\right] _{h^{\ast }F}\in
\Gamma \left( E,\pi ,M\right) ,~$for any $S,T\in \Gamma \left( E,\pi
,M\right) .$

\textbf{Proposition 4.3 }\emph{If }$\left( E,\pi ,M\right) $\emph{\ is an
IDS of the generalized Lie algebroid }%
\begin{equation*}
\left( \left( F,\nu ,N\right) ,\left[ ,\right] _{F,h},\left( \rho ,\eta
\right) \right)
\end{equation*}%
\emph{and }$\left\{ S_{1},...,S_{r}\right\} $\emph{\ is a base for the }$%
\mathcal{F}\left( M\right) $\emph{-submodule }$\left( \Gamma \left( E,\pi
,M\right) ,+,\cdot \right) $\emph{\ then }$\left( E,\pi ,M\right) $\emph{\
is involutive if and only if }$\left[ S_{a},S_{b}\right] _{h^{\ast }F}\in
\Gamma \left( E,\pi ,M\right) ,~$for any $a,b\in \overline{1,r}.$

\textbf{Theorem 4.2 (}of Frobenius type)\textbf{\ }\emph{Let }$\left( E,\pi
,M\right) $\emph{\ be an IDS of the generalized Lie algebroid }$\left(
\left( F,\nu ,N\right) ,\left[ ,\right] _{F,h},\left( \rho ,\eta \right)
\right) .$\emph{\ If }$\left\{ \Theta ^{r+1},...,\Theta ^{p}\right\} $\emph{%
\ is a base for the }$\mathcal{F}\left( M\right) $\emph{-submodule }$\left(
\Gamma \left( E^{0},\pi ^{0},M\right) ,+,\cdot \right) $\emph{, then the IDS 
}$\left( E,\pi ,M\right) $\emph{\ is involutive if and only if it exists }%
\begin{equation*}
\Omega _{\beta }^{\alpha }\in \Lambda ^{1}\left( h^{\ast }F,h^{\ast }\nu
,M\right) ,~\alpha ,\beta \in \overline{r+1,p}
\end{equation*}%
\emph{such that} 
\begin{equation*}
d^{h^{\ast }F}\Theta ^{\alpha }=\Sigma _{\beta \in \overline{r+1,p}}\Omega
_{\beta }^{\alpha }\wedge \Theta ^{\beta },~\alpha \in \overline{r+1,p}.
\end{equation*}

\emph{Proof: }Let $\left\{ S_{1},...,S_{r}\right\} $ be a base for the $%
\mathcal{F}\left( M\right) $-submodule $\left( \Gamma \left( E,\pi ,M\right)
,+,\cdot \right) $

Let $\left\{ S_{r+1},...,S_{p}\right\} \in \Gamma \left( h^{\ast }F,h^{\ast
}\nu ,M\right) $ such that 
\begin{equation*}
\left\{ S_{1},...,S_{r},S_{r+1},...,S_{p}\right\}
\end{equation*}%
is a base for the $\mathcal{F}\left( M\right) $-module 
\begin{equation*}
\left( \Gamma \left( h^{\ast }F,h^{\ast }\nu ,M\right) ,+,\cdot \right) .
\end{equation*}

Let $\Theta ^{1},...,\Theta ^{r}\in \Gamma \left( \overset{\ast }{h^{\ast }F}%
,\overset{\ast }{h^{\ast }\nu },M\right) $ such that 
\begin{equation*}
\left\{ \Theta ^{1},...,\Theta ^{r},\Theta ^{r+1},...,\Theta ^{p}\right\}
\end{equation*}
is a base for the $\mathcal{F}\left( M\right) $-module 
\begin{equation*}
\left( \Gamma \left( \overset{\ast }{h^{\ast }F},\overset{\ast }{h^{\ast
}\nu },M\right) ,+,\cdot \right) .
\end{equation*}

For any $a,b\in \overline{1,r}$ and $\alpha ,\beta \in \overline{r+1,p}$, we
have the equalities:%
\begin{equation*}
\begin{array}{ccc}
\Theta ^{a}\left( S_{b}\right) & = & \delta _{b}^{a} \\ 
\Theta ^{a}\left( S_{\beta }\right) & = & 0 \\ 
\Theta ^{\alpha }\left( S_{b}\right) & = & 0 \\ 
\Theta ^{\alpha }\left( S_{\beta }\right) & = & \delta _{\beta }^{\alpha }%
\end{array}%
\end{equation*}

We remark that the set of the $2$-forms%
\begin{equation*}
\left\{ \Theta ^{a}\wedge \Theta ^{b},\Theta ^{a}\wedge \Theta ^{\beta
},\Theta ^{\alpha }\wedge \Theta ^{\beta },~a,b\in \overline{1,r}\wedge
\alpha ,\beta \in \overline{r+1,p}\right\}
\end{equation*}%
is a base for the $\mathcal{F}\left( M\right) $-module 
\begin{equation*}
\left( \Lambda ^{2}\left( h^{\ast }F,h^{\ast }\nu ,M\right) ,+,\cdot \right)
.
\end{equation*}

Therefore, we have%
\begin{equation*}
d^{h^{\ast }F}\Theta ^{\alpha }=\Sigma _{b<c}A_{bc}^{\alpha }\Theta
^{b}\wedge \Theta ^{c}+\Sigma _{b,\gamma }B_{b\gamma }^{\alpha }\Theta
^{b}\wedge \Theta ^{\gamma }+\Sigma _{\beta <\gamma }C_{\beta \gamma
}^{\alpha }\Theta ^{\beta }\wedge \Theta ^{\gamma },\leqno\left( 1\right)
\end{equation*}%
where, $A_{bc}^{\alpha },B_{b\gamma }^{\alpha }$ and $C_{\beta \gamma
}^{\alpha },~a,b,c\in \overline{1,r}\wedge \alpha ,\beta ,\gamma \in 
\overline{r+1,p}$ are real local functions such that $A_{bc}^{\alpha
}=-A_{cb}^{\alpha }$ and $C_{\beta \gamma }^{\alpha }=-C_{\gamma \beta
}^{\alpha }.$

Using the formula%
\begin{equation*}
d^{h^{\ast }F}\Theta ^{\alpha }\left( S_{b},S_{c}\right) =\Gamma \left( 
\overset{h^{\ast }F}{\rho },Id_{M}\right) S_{b}\left( \Theta ^{\alpha
}\left( S_{c}\right) \right) -\Gamma \left( \overset{h^{\ast }F}{\rho }%
,Id_{M}\right) S_{c}\left( \Theta ^{\alpha }\left( S_{b}\right) \right)
-\Theta ^{\alpha }\left( \left[ S_{b},S_{c}\right] _{h^{\ast }F}\right) ,%
\leqno\left( 2\right)
\end{equation*}%
we obtain that 
\begin{equation*}
A_{bc}^{\alpha }=-\Theta ^{\alpha }\left( \left[ S_{b},S_{c}\right]
_{h^{\ast }F}\right) ,~\forall \left( b,c\in \overline{1,r}\wedge \alpha \in 
\overline{r+1,p}\right) .\leqno\left( 3\right)
\end{equation*}

We admit that $\left( E,\pi ,M\right) $ is an involutive \emph{IDS} of the
generalized Lie algebroid $\left( \left( F,\nu ,N\right) ,\left[ ,\right]
_{F,h},\left( \rho ,\eta \right) \right) .$

As 
\begin{equation*}
\left[ S_{b},S_{c}\right] _{h^{\ast }F}\in \Gamma \left( E,\pi ,M\right)
,~\forall b,c\in \overline{1,r}
\end{equation*}%
it results that 
\begin{equation*}
\Theta ^{\alpha }\left( \left[ S_{b},S_{c}\right] _{h^{\ast }F}\right)
=0,~\forall \left( b,c\in \overline{1,r}\wedge \alpha \in \overline{r+1,p}%
\right) .
\end{equation*}

Therefore, 
\begin{equation*}
A_{bc}^{\alpha }=0,~\forall \left( b,c\in \overline{1,r}\wedge \alpha \in 
\overline{r+1,p}\right)
\end{equation*}%
and we obtain%
\begin{equation*}
\begin{array}{ccl}
d^{h^{\ast }F}\Theta ^{\alpha } & = & \Sigma _{b,\gamma }B_{b\gamma
}^{\alpha }\Theta ^{b}\wedge \Theta ^{\gamma }+\frac{1}{2}C_{\beta \gamma
}^{\alpha }\Theta ^{\beta }\wedge \Theta ^{\gamma } \\ 
& = & \left( B_{b\gamma }^{\alpha }\Theta ^{b}+\frac{1}{2}C_{\beta \gamma
}^{\alpha }\Theta ^{\beta }\right) \wedge \Theta ^{\gamma }.%
\end{array}%
\end{equation*}

As 
\begin{equation*}
\Omega _{\gamma }^{\alpha }\overset{put}{=}B_{b\gamma }^{\alpha }\Theta ^{b}+%
\frac{1}{2}C_{\beta \gamma }^{\alpha }\Theta ^{\beta }\in \Lambda ^{1}\left(
h^{\ast }F,h^{\ast }\nu ,M\right) ,~\forall \alpha ,\beta \in \overline{r+1,p%
}
\end{equation*}%
it results the first implication.

Conversely, we admit that it exists 
\begin{equation*}
\Omega _{\beta }^{\alpha }\in \Lambda ^{1}\left( h^{\ast }F,h^{\ast }\nu
,M\right) ,~\alpha ,\beta \in \overline{r+1,p}
\end{equation*}%
such that 
\begin{equation*}
d^{h^{\ast }F}\Theta ^{\alpha }=\Sigma _{\beta \in \overline{r+1,p}}\Omega
_{\beta }^{\alpha }\wedge \Theta ^{\beta },~\forall \alpha \in \overline{%
r+1,p}.\leqno\left( 4\right)
\end{equation*}

Using the affirmations $\left( 1\right) ,\left( 2\right) $ and $\left(
4\right) $ we obtain that 
\begin{equation*}
A_{bc}^{\alpha }=0,~\forall \left( b,c\in \overline{1,r}\wedge \alpha \in 
\overline{r+1,p}\right) .
\end{equation*}

Using the affirmation $\left( 3\right) $, we obtain 
\begin{equation*}
\Theta ^{\alpha }\left( \left[ S_{b},S_{c}\right] _{h^{\ast }F}\right)
=0,~\forall \left( b,c\in \overline{1,r}\wedge \alpha \in \overline{r+1,p}%
\right) .
\end{equation*}

Therefore, 
\begin{equation*}
\left[ S_{b},S_{c}\right] _{h^{\ast }F}\in \Gamma \left( E,\pi ,M\right)
,~\forall b,c\in \overline{1,r}.
\end{equation*}

Using the \emph{Proposition 4.2}, we obtain the second implication.\hfill 
\emph{q.e.d.}\medskip

Let $\left( \left( h^{\ast }F,h^{\ast }\nu ,M\right) ,\left[ ,\right]
_{h^{\ast }F},\left( \overset{h^{\ast }F}{\rho },Id_{M}\right) \right) $ be
the pull-back Lie algebroid of the generalized Lie algebroid $\left( \left(
F,\nu ,N\right) ,\left[ ,\right] _{F,h},\left( \rho ,\eta \right) \right) $.

\textbf{Definition 4.3 }Any ideal $\left( \mathcal{I},+,\cdot \right) $ of
the exterior differential algebra of the pull-back Lie algebroid $\left(
\left( h^{\ast }F,h^{\ast }\nu ,M\right) ,\left[ ,\right] _{h^{\ast
}F},\left( \overset{h^{\ast }F}{\rho },Id_{M}\right) \right) $ closed under
differentiation operator $d^{h^{\ast }F},$ namely $d^{h^{\ast }F}\mathcal{%
I\subseteq I},$ will be called \emph{differential ideal of the generalized
Lie algebroid }$\left( \left( F,\nu ,N\right) ,\left[ ,\right] _{F,h},\left(
\rho ,\eta \right) \right) .$

In particular, if $h=Id_{N}=\eta $, then we obtain the definition of the
differential ideal of a Lie algebroid.(see$\left[ 2\right] $)

\textbf{Definition 4.5 }Let $\left( \mathcal{I},+,\cdot \right) $ be a
differential ideal of the generalized Lie algebroid $\left( \left( F,\nu
,N\right) ,\left[ ,\right] _{F,h},\left( \rho ,\eta \right) \right) $.

If it exists an \emph{IDS }$\left( E,\pi ,M\right) $ such that for all $k\in 
\mathbb{N}^{\ast }$ and $\omega \in \mathcal{I}\cap \Lambda ^{k}\left(
h^{\ast }F,h^{\ast }\nu ,M\right) $ we have $\omega \left(
u_{1},...,u_{k}\right) =0,$ for any $u_{1},...,u_{k}\in \Gamma \left( E,\pi
,M\right) ,$ then we will say that $\left( \mathcal{I},+,\cdot \right) $%
\emph{\ is an exterior differential system (EDS) of the generalized Lie
algebroid }%
\begin{equation*}
\left( \left( F,\nu ,N\right) ,\left[ ,\right] _{F,h},\left( \rho ,\eta
\right) \right) .
\end{equation*}

In particular, if $h=Id_{N}=\eta $, then we obtain the definition of the 
\emph{EDS }of a Lie algebroid.(see$\left[ 2\right] $)

\textbf{Theorem 4.3 (}of Cartan type) \emph{The IDS }$\left( E,\pi ,M\right) 
$\emph{\ of the generalized Lie algebroid }$\left( \left( F,\nu ,N\right) ,%
\left[ ,\right] _{F,h},\left( \rho ,\eta \right) \right) $\emph{\ is
involutive, if and only if the ideal generated by the }$\mathcal{F}\left(
M\right) $\emph{-submodule }$\left( \Gamma \left( E^{0},\pi ^{0},M\right)
,+,\cdot \right) $\emph{\ is an EDS of the same generalized Lie algebroid.}

\emph{Proof. }Let $\left( E,\pi ,M\right) $ be an involutive \emph{IDS} of
the generalized Lie algebroid 
\begin{equation*}
\left( \left( F,\nu ,N\right) ,\left[ ,\right] _{F,h},\left( \rho ,\eta
\right) \right) .
\end{equation*}

Let $\left\{ \Theta ^{r+1},...,\Theta ^{p}\right\} $ be a base for the $%
\mathcal{F}\left( M\right) $-submodule $\left( \Gamma \left( E^{0},\pi
^{0},M\right) ,+,\cdot \right) .$

We know that 
\begin{equation*}
\mathcal{I}\left( \Gamma \left( E^{0},\pi ^{0},M\right) \right) =\cup _{q\in 
\mathbb{N}}\left\{ \Omega _{\alpha }\wedge \Theta ^{\alpha },~\left\{ \Omega
_{r+1},...,\Omega _{p}\right\} \subset \Lambda ^{q}\left( h^{\ast }F,h^{\ast
}\nu ,M\right) \right\} .
\end{equation*}

Let $q\in \mathbb{N}$ and $\left\{ \Omega _{r+1},...,\Omega _{p}\right\}
\subset \Lambda ^{q}\left( h^{\ast }F,h^{\ast }\nu ,M\right) $ be arbitrary.

Using the \emph{Theorems 3.8 and 3.10} we obtain 
\begin{equation*}
\begin{array}{ccl}
d^{h^{\ast }F}\left( \Omega _{\alpha }\wedge \Theta ^{\alpha }\right) & = & 
d^{h^{\ast }F}\Omega _{\alpha }\wedge \Theta ^{\alpha }+\left( -1\right)
^{q+1}\Omega _{\beta }\wedge d^{h^{\ast }F}\Theta ^{\beta } \\ 
& = & \left( d^{h^{\ast }F}\Omega _{\alpha }+\left( -1\right) ^{q+1}\Omega
_{\beta }\wedge \Omega _{\alpha }^{\beta }\right) \wedge \Theta ^{\alpha }.%
\end{array}%
\end{equation*}

As 
\begin{equation*}
d^{h^{\ast }F}\Omega _{\alpha }+\left( -1\right) ^{q+1}\Omega _{\beta
}\wedge \Omega _{\alpha }^{\beta }\in \Lambda ^{q+2}\left( h^{\ast
}F,h^{\ast }\nu ,M\right)
\end{equation*}%
it results that 
\begin{equation*}
d^{h^{\ast }F}\left( \Omega _{\beta }\wedge \Theta ^{\beta }\right) \in 
\mathcal{I}\left( \Gamma \left( E^{0},\pi ^{0},M\right) \right)
\end{equation*}

Therefore, 
\begin{equation*}
d^{h^{\ast }F}\mathcal{I}\left( \Gamma \left( E^{0},\pi ^{0},M\right)
\right) \subseteq \mathcal{I}\left( \Gamma \left( E^{0},\pi ^{0},M\right)
\right) .
\end{equation*}

Conversely, let $\left( E,\pi ,M\right) $ be an \emph{IDS} of the
generalized Lie algebroid 
\begin{equation*}
\left( \left( F,\nu ,N\right) ,\left[ ,\right] _{F,h},\left( \rho ,\eta
\right) \right)
\end{equation*}%
such that the $\mathcal{F}\left( M\right) $-submodule $\left( \mathcal{I}%
\left( \Gamma \left( E^{0},\pi ^{0},M\right) \right) ,+,\cdot \right) $ is
an \emph{EDS} of the generalized Lie algebroid $\left( \left( F,\nu
,N\right) ,\left[ ,\right] _{F,h},\left( \rho ,\eta \right) \right) .$

Let $\left\{ \Theta ^{r+1},...,\Theta ^{p}\right\} $ be a base for the $%
\mathcal{F}\left( M\right) $-submodule $\left( \Gamma \left( E^{0},\pi
^{0},M\right) ,+,\cdot \right) .$ As 
\begin{equation*}
d^{h^{\ast }F}\mathcal{I}\left( \Gamma \left( E^{0},\pi ^{0},M\right)
\right) \subseteq \mathcal{I}\left( \Gamma \left( E^{0},\pi ^{0},M\right)
\right)
\end{equation*}%
it results that it exists 
\begin{equation*}
\Omega _{\beta }^{\alpha }\in \Lambda ^{1}\left( h^{\ast }F,h^{\ast }\nu
,M\right) ,~\alpha ,\beta \in \overline{r+1,p}
\end{equation*}%
such that 
\begin{equation*}
d^{h^{\ast }F}\Theta ^{\alpha }=\Sigma _{\beta \in \overline{r+1,p}}\Omega
_{\beta }^{\alpha }\wedge \Theta ^{\beta }\in \mathcal{I}\left( \Gamma
\left( E^{0},\pi ^{0},M\right) \right) .
\end{equation*}

Using the \emph{Theorem 4.2}, it results that $\left( E,\pi ,M\right) $ is
an involutive \emph{IDS.}\hfill \emph{q.e.d.}\medskip

\bigskip

\section{A new direction by research}

\ \ \ \ 

We know that the set of morphisms of 
\begin{equation*}
\left( \left( F,\nu ,N\right) ,\left[ ,\right] _{F,h},\left( \rho ,\eta
\right) \right)
\end{equation*}%
source and 
\begin{equation*}
\left( \left( F^{\prime },\nu ^{\prime },N^{\prime }\right) ,\left[ ,\right]
_{F^{\prime },h^{\prime }},\left( \rho ^{\prime },\eta ^{\prime }\right)
\right)
\end{equation*}%
target is the set 
\begin{equation*}
\begin{array}{c}
\left\{ \left( \varphi ,\varphi _{0}\right) \in \mathbf{B}^{\mathbf{v}%
}\left( \left( F,\nu ,N\right) ,\left( F^{\prime },\nu ^{\prime },N^{\prime
}\right) \right) \right\}%
\end{array}%
\end{equation*}%
such that $\varphi _{0}\in Iso_{\mathbf{Man}}\left( N,N^{\prime }\right) $
and the $\mathbf{Mod}$-morphism $\Gamma \left( \varphi ,\varphi _{0}\right) $
is a $\mathbf{LieAlg}$-morphism of 
\begin{equation*}
\left( \Gamma \left( F,\nu ,N\right) ,+,\cdot ,\left[ ,\right] _{F,h}\right)
\end{equation*}%
source and 
\begin{equation*}
\left( \Gamma \left( F^{\prime },\nu ^{\prime },N^{\prime }\right) ,+,\cdot
, \left[ ,\right] _{F^{\prime },h^{\prime }}\right)
\end{equation*}%
target. (see: $\left[ 1\right] $)

We can define \emph{the simplectic space }as being a pair 
\begin{equation*}
\begin{array}{c}
\left( \left( \left( F,\nu ,N\right) ,\left[ ,\right] _{F,h},\left( \rho
,\eta \right) \right) ,\omega \right)%
\end{array}%
\end{equation*}%
consisting of a generalized Lie algebroid $\left( \left( F,\nu ,N\right) ,%
\left[ ,\right] _{F,h},\left( \rho ,\eta \right) \right) $ and a
nondegenerate close $2$-form $\omega \in \Lambda ^{2}\left( F,\nu ,N\right)
. $

If 
\begin{equation*}
\begin{array}{c}
\left( \left( \left( F^{\prime },\nu ^{\prime },N^{\prime }\right) ,\left[ ,%
\right] _{F^{\prime },h^{\prime }},\left( \rho ^{\prime },\eta ^{\prime
}\right) \right) ,\omega ^{\prime }\right)%
\end{array}%
\end{equation*}%
is an another simplectic space, then we can define the set of morphisms of 
\begin{equation*}
\begin{array}{c}
\left( \left( \left( F,\nu ,N\right) ,\left[ ,\right] _{F,h},\left( \rho
,\eta \right) \right) ,\omega \right)%
\end{array}%
\end{equation*}%
source and 
\begin{equation*}
\begin{array}{c}
\left( \left( \left( F^{\prime },\nu ^{\prime },N^{\prime }\right) ,\left[ ,%
\right] _{F^{\prime },h^{\prime }},\left( \rho ^{\prime },\eta ^{\prime
}\right) \right) ,\omega ^{\prime }\right)%
\end{array}%
\end{equation*}%
target as being the set 
\begin{equation*}
\begin{array}{c}
\left\{ \left( \varphi ,\varphi _{0}\right) \in \mathbf{GLA}\left( \left(
\left( F,\nu ,N\right) ,\left[ ,\right] _{F,h},\left( \rho ,\eta \right)
\right) ,\left( \left( F^{\prime },\nu ^{\prime },N^{\prime }\right) ,\left[
,\right] _{F^{\prime },h^{\prime }},\left( \rho ^{\prime },\eta ^{\prime
}\right) \right) \right) \right\}%
\end{array}%
\end{equation*}%
such that $\left( \varphi ,\varphi _{0}\right) ^{\ast }\left( \omega
^{\prime }\right) =\omega .$

So, we can discuss about the category of simplectic spaces as being a
subcategory of the category of generalized Lie algebroids. The study of the
geometry of objects of this category is a new direction by research.

Very interesting will be a result of Darboux type in this general framework
and the connections with the Poisson bracket.

\bigskip

\section{Acknowledgment}

\addcontentsline{toc}{section}{Acknowledgment}

I would like to thank R\u{a}dine\c{s}ti-Gorj Cultural Scientifique Society
for financial support. In memory of Prof. Dr. Gheorghe RADU and Acad. Dr.
Doc. Cornelius RADU. Dedicated to Prof. Dr. Sorin Vasile SABAU from Tokai
University, Japan.

\hfill


\bigskip \addcontentsline{toc}{section}{References}

\begin{thebibliography}{99}
\bibitem{1} C. M. Arcus, \emph{Generalized Lie algebroids and connections
over pair of diffeomorphic manifolds}, Journal of Generalized Lie Theory and
Applications, vol. 7, (2012), Article ID G111202, 32 pages,
doi:10.4303/glta/G111202.

\bibitem{2} C. M. Arcus, \emph{Interior and exterior differential systems
for Lie algebroids}, Advances in Pure Mathematics, Vol. I, No. 5, (2011),
245-249, doi:10.4236/apm.2011.15004.

\bibitem{3} R.L. Bryant, S.S. Chern, R.B. Gardner, H.L. Goldschmidt, and
P.A. Griffiths, \emph{Exterior Differential Systems,} Springer-Verlag, 1991.

\bibitem{4} J. Grabowski, P. Urbanski, \emph{Lie algebroids and
Poisson-Nijenhuis structures,} Rep. Math. Phys. \textbf{40} (1997), 196-208,
doi:10.1016/S0034-4877(97).85916-2.

\bibitem{5} J. Grabowski, P. Urbanski, \emph{Algebroids-General differential
calculi on vector bundles,} Geom. Phys. \textbf{31} (1999), 111-141.

\bibitem{6} P Griffiths, \emph{Exterior Differential Systems and the
Calculus of Variations,} Progr. Math., No. 25, Birkh\"{a}user, Boston, MA,
1983.

\bibitem{7} T. A. Ivey and J. M. Landsberg, \emph{Cartan for Beginners:
Differential Geometry via Moving Frames and Exterior Differential Systems,}
Graduate Texts in Mathematics, American Mathematical Society, 2003.

\bibitem{8} N. Kamran, \emph{An elementary introduction to exterior
differential systems,} In \textquotedblleft Geometric approaches to
differential equations (Canberra, 1995)\textquotedblright , volume 15 of
Austral. Math. Soc. Lect. Ser., Cambridge Univ. Press, Cambridge, (2000),
100--115.

\bibitem{9} M. de Leon, \emph{Methods on Differential Geometry in Analitical
Mechanics,} North-Holand, Amsterdam, 1989.

\bibitem{10} C. M. Marle, \emph{Lie Algebroids and Lie Pseudoalgebras,}
Mathematics \& Physical Sciences, Vol. 27, No. 2, (2008), 97-147. \ \ 

\bibitem{11} C. M. Marle, \emph{Calculus on Lie algebroids, Lie groupoids
and Poisson manifolds,} arXiv:0806.0919v3 [math.DG] (2008)

\bibitem{12} L. I. Nicolescu, \emph{Lectures on Geometry of Manifolds, }%
World Sciantific, Singapore, 1996, doi:10.1142/9789814261012.
\end{thebibliography}
\end{document}